\documentclass[a4paper,fleqn]{cas-sc}

\usepackage[authoryear]{natbib}
\usepackage[figuresright]{rotating}
\usepackage{subfigure}
\usepackage{xcolor} 
\usepackage{colortbl}  

\begin{document}
\let\WriteBookmarks\relax
\def\floatpagepagefraction{1}
\def\textpagefraction{.001}
\let\printorcid\relax   


\shorttitle{Collaborative Charging Scheduling via B3M}

\shortauthors{Zhou et~al.}

\title [mode = title]{Collaborative Charging Scheduling via Balanced Bounding Box Methods}

\author[1]{Fangting Zhou$^*$}

\author[2]{Balázs Kulcsár}

\author[1]{Jiaming Wu}

\affiliation[1]{organization={Architecture and Civil Engineering, Chalmers University of Technology},
    city={Gothenburg},
    country={Sweden}}

\affiliation[2]{organization={Electrical Engineering, Chalmers University of Technology},
    city={Gothenburg},
    country={Sweden}}

\begin{abstract}
Electric mobility faces several challenges, most notably the high cost of infrastructure development and the underutilization of charging stations. The concept of shared charging offers a promising solution. The paper explores sustainable urban logistics through horizontal collaboration between two fleet operators and addresses a scheduling problem for the shared use of charging stations. To tackle this, the study formulates a collaborative scheduling problem as a bi-objective nonlinear integer programming model, in which each company aims to minimize its own costs, creating inherent conflicts that require trade-offs. The Balanced Bounding Box Methods (B3Ms) are introduced in order to efficiently derive the efficient frontier, identifying a reduced set of representative solutions. These methods enhance computational efficiency by selectively disregarding closely positioned and competing solutions, preserving the diversity and representativeness of the solutions over the efficient frontier. To determine the final solution and ensure balanced collaboration, cooperative bargaining methods are applied. Numerical case studies demonstrate the viability and scalability of the developed methods, showing that the B3Ms can significantly reduce computational time while maintaining the integrity of the frontier. These methods, along with cooperative bargaining, provide an effective framework for solving various bi-objective optimization problems, extending beyond the collaborative scheduling problem presented here.
\end{abstract}

\begin{keywords}
Urban logistics \sep Collaborative scheduling \sep Electric vehicle \sep Shared charging \sep Balanced Bounding Box Method \sep Nash bargaining
\end{keywords}

\maketitle

\renewcommand{\thefootnote}{}%
\footnotetext{* Corresponding author. Email: fangting@chalmers.se}%
\renewcommand{\thefootnote}{\arabic{footnote}}  

\section{Introduction}


Electric vehicles (EVs) are widely recognized as a key solution for decarbonizing transportation, offering up to 50\% lower greenhouse gas emissions and reduced energy consumption compared to gasoline vehicles \citep{wang2019, nealer2015, kumar2021}. As a result, global EV adoption is accelerating rapidly. The global EV stock is projected to surpass 500 million by 2035, growing at an average rate of 23\% annually \citep{IEA2024}. This surge places unprecedented pressure on charging infrastructure, especially public charging networks, which are expected to expand significantly to support this growth \citep{bauer2021}. Yet, ensuring that these networks are efficiently coordinated and economically viable remains a persistent challenge.


However, expanding public charging infrastructure remains a costly and complex endeavor, involving substantial investment in installation, maintenance, and grid upgrades \citep{ACEA2024}. Public charging prices are significantly higher than private ones, thereby imposing a greater financial burden on users, particularly fleet operators \citep{kampshoff2022}. Moreover, low utilization rates further undermine economic viability. In 2022, public fast-charging stations in the U.S. averaged only 7.5\% utilization, far below the 15\% threshold typically required for financial feasibility \citep{McKinsey2022}. 
Moreover, recent studies highlight that charging infrastructure expansion is often constrained by scheduling feasibility, user behavior, and grid limitations, suggesting that physical development alone may be insufficient \citep{li2024, rizopoulos2025}. These challenges underscore the need for smarter scheduling strategies to improve utilization and cost-effectiveness.


In response to these challenges, shared charging has emerged as a promising approach to improve infrastructure efficiency and reduce investment costs by enabling multiple users to access the same facilities. Recent studies show that such models not only increase station utilization but also generate additional revenue for infrastructure providers \citep{gao2021, hu2021, melander2023, cai2024}. Platforms like Share\&Charge and Co Charger have adopted this concept, allowing individuals and businesses to rent out underutilized chargers via digital marketplaces \citep{vanrykel2018,CoCharger2024}.


Existing research on shared charging primarily focuses on pricing strategies and market mechanisms for revenue generation. For example, \citet{gao2021} propose an auction-based pricing model, while \citet{cai2024} develop a ride-sourcing charger-sharing program. Moreover, \citet{ji2023} proposed a time-sharing model that allows electric buses and private vehicles to jointly use charging stations, improving utilization and cost-efficiency. However, their model assumes a dominant operator and does not address fairness or coordination among multiple self-interested users.
More broadly, existing studies emphasize economic incentives without fully addressing operational issues such as user competition, scheduling conflicts, or cost allocation.
Meanwhile, fleet-based EV scheduling has been studied extensively \citep{kang2013, zakariazadeh2014, wu2020}, some using multi-objective models \citep{wu2020, hajforoosh2015}, but mainly for individual users rather than collaborative use of shared infrastructure.


Despite their potential, shared charging systems pose operational challenges due to multi-user interactions and limited scheduling resources. Most models optimize single-user objectives and lack cooperative mechanisms to allocate access and cost fairly. As users compete for chargers, minimizing conflicts while ensuring equitable outcomes remains difficult. Moreover, existing coordination frameworks often fall short in balancing efficiency and fairness, limiting scalability and real-world applicability. Addressing these shortcomings is crucial to unlocking the full potential of shared infrastructure.

To bridge this gap, this study proposes a ``rent and share'' scheduling approach, where EV fleet operators collaboratively rent and share charging infrastructure. By leveraging existing stations through coordinated scheduling, this approach reduces costs, enhances resource utilization, and minimizes operational conflicts.
To facilitate fair and efficient collaboration between companies, we develop a collaborative scheduling framework based on horizontal collaboration. This concept allows companies operating at the same level in the market to share resources and coordinate activities to improve infrastructure efficiency and operational performance \citep{cruijssen2007}. In this study, two EV fleet operators rent charging stations from third-party providers instead of maintaining dedicated infrastructure, ensuring balanced scheduling efficiency and cost allocation.

Building on this framework, we introduce a bi-objective optimization model integrated with game theory, jointly optimizing two interconnected decisions: (i) determining which charging stations each company should rent, and (ii) coordinating the scheduling of charging activities across the shared network. 
These decisions are jointly optimized within the framework that balances cost efficiency and fairness, ensuring equitable resource allocation while maximizing overall operational efficiency. To efficiently solve the bi-objective optimization model, we develop the Balanced Bounding Box Method (B3M), which delineates the non-dominated frontier by filtering out closely positioned solutions to generate a compact, representative subset, thereby improving computational efficiency and solution diversity. Additionally, the framework integrates cooperative bargaining mechanisms to model strategic interactions between companies, guiding them toward mutually beneficial agreements that minimize costs while effectively managing trade-offs in their competing objectives.
Numerical case studies validate the scalability and effectiveness of our approach, demonstrating significant computational efficiency while maintaining solution quality.

The main contributions of this paper are as follows:

\begin{itemize}
    \setlength{\itemsep}{0pt}
    \setlength{\parsep}{0pt}
    \setlength{\parskip}{0pt}
    \item A collaborative scheduling problem is defined, which coordinates the charging needs between two companies that rent and share charging stations.
    \item A bi-objective integer programming model is proposed, where each company minimizes its own costs, leading to inherent trade-offs.
    \item The B3Ms are developed to efficiently derive an efficient frontier with a reduced set of representative non-dominated solutions.
    \item Cooperative bargaining approaches are integrated to select a final solution from the non-dominated set, ensuring fair and robust decision-making for both companies.
\end{itemize}

The remainder of this paper is organized as follows. Section \ref{Section1} reviews the relevant literature. Section \ref{Section2} outlines the problem and formulates the mathematical programming model. Section \ref{Section3} develops B3Ms to identify the non-dominated frontier, and presents two cooperative bargaining methods for determining the final agreement point. Section \ref{Section4} presents the experimental study and numerical results. Finally, Section \ref{Section5} concludes the paper and suggests directions for future research.

\section{Literature review} \label{Section1}

Developing EV charging infrastructure faces several challenges, including planning \citep{mak2013}, the impact of public charging networks \citep{levinson2018}, and the allocation of investment costs \citep{kumar2021}. Sharing these costs across stakeholders with vested interests offers a potential solution to these issues \citep{bauer2021}. However, while cost-sharing is promising, optimizing the usage and coordination of shared charging facilities remains an open question, especially in complex urban environments where demand is high and resources are limited.

\subsection*{Shared charging stations}

Research on shared charging stations remains limited but highlights several innovative approaches for optimizing infrastructure use and market dynamics. \cite{gao2021} propose a price-based iterative double auction model for charger-sharing markets, where charger owners rent underutilized chargers to EV drivers, creating a two-sided market regulated by auction. This model accounts for factors such as charger location, availability, service costs, and driver preferences, aiming to maximize social welfare and ensure mutual benefits for both charger owners and EV drivers. Similarly, \cite{cai2024} introduce a shared charging program where a ride-sourcing platform acts as an intermediary, renting charging posts from private owners to serve ride-sourcing drivers. This utility-based model analyzes interactions among post owners, drivers, and riders, exploring optimal pricing and decision-making strategies to improve social welfare and market outcomes in both the charging and ride-sourcing sectors.
In addition to these pricing-focused approaches, \citet{ji2023} develop an optimization model that integrates electric bus (EB) scheduling with charging facility sharing, allowing electric cars (ECs) to access EB chargers during specified time windows. Their model jointly minimizes EB operational costs and EC waiting time while maximizing charging revenue, offering a coordinated time-sharing strategy to improve utilization. However, this model assumes a single dominant operator and does not address interactions or conflicts between multiple users.

Despite growing interest in shared charging, existing research remains largely limited to pricing and market structures, with little attention given to scheduling mechanisms that can ensure fairness, resolve conflicts, and support coordination among multiple competing users. This gap highlights the need to draw on scheduling-based approaches to enable more efficient and equitable shared infrastructure use.

\subsection*{EV charging scheduling}

EV charging scheduling is a critical and complex issue that has drawn significant research interest \citep{mukherjee2014}. Researchers have developed scheduling schemes for various charging stations, including standard EV stations \citep{kang2013, zakariazadeh2014, zhang2018optimal, he2020, das2020}, bus networks \citep{duan2023, zhou2024electric}, battery swapping stations \citep{sarker2014, you2015, wu2017, tan2023}, and stations integrated with renewable energy sources \citep{kabir2020, long2021, das2021}.

More specifically, \cite{wei2017} focus on optimizing EV charging schedules within parking garages to maximize utility for both operators and customers, incorporating time-of-use pricing. They propose a multi-charging system that accounts for real-world battery behaviors and uses intelligent management strategies. Similarly, \cite{wu2020} examine a time-slotted model aimed at minimizing total charging costs and reducing peak load. Building on this, \cite{saner2022} introduce a hierarchical multi-agent system to coordinate multiple EV charging stations. In this system, higher-level agents generate optimized control signals for lower-level agents, enabling effective coordination without direct communication or prior knowledge of EV arrivals.

Other studies focus on scheduling specific EV fleets. For example, \cite{bao2023} address scheduling for an electric bus fleet, determining optimal charging within defined time windows. Each window includes charging slots of equal duration, aiming to minimize total fleet costs while considering consumption rates, tariffs, SOC limits, and station capacities.

Additionally, \cite{wu2023} propose a centralized EV charging scheduling strategy under time-of-use pricing, aiming to minimize charging costs while adhering to constraints related to charging pile availability and station power capacity. This model is designed to optimize scheduling for a shared charging environment, rather than a specific fleet, by adapting to user demand and managing instantaneous power constraints.

Beyond scheduling, research has also focused on the assignment of EVs to charging stations to optimize infrastructure utilization. For example, \cite{chen2013} present a model for EV charging station location assignment that minimizes access costs for users while penalizing unmet demand. Addressing the need for dynamic, real-time solutions, \cite{ma2021} propose an online model that dynamically assigns vehicles to charging stations within ridesharing contexts. This model focuses on optimizing the placement of fast charging stations to ensure efficient charging management. By tackling various aspects of EV-to-station assignment, these studies offer insights into improving charging infrastructure utilization and efficiency.

To further balance the diverse objectives in EV charging management, some studies use multi-objective optimization to handle conflicting interests. Most commonly, these studies address two conflicting objectives, referred to as bi-objective optimization, such as minimizing costs \citep{zakariazadeh2014, hajforoosh2015, wu2020}, emissions \citep{zakariazadeh2014}, and peak load \citep{wu2020}, or maximizing delivered power \citep{hajforoosh2015} and battery stock levels \citep{wu2017}. More complex models consider three objectives, including electricity cost, battery degradation cost, and grid net exchange \citep{das2020, das2021}. Coordinating among stakeholders with such conflicting objectives remains challenging. In this paper, we focus on horizontal collaboration, involving EV manufacturers or companies with EV fleets that need to coordinate their charging schedules.

Despite extensive research on optimizing EV charging schedules and station assignments, collaborative scheduling mechanisms have received limited attention, especially for horizontal collaboration among stakeholders such as EV manufacturers and fleet operators with competing objectives. This highlights the necessity for coordination strategies to improve decision-making and resource sharing across entities.

\subsection*{Solution methods}

Building on the concept of horizontal collaboration among EV fleet operators, it is essential to address the multi-objective nature of charging schedules. Multi-objective optimization balances conflicting goals, enabling high-quality solutions in real-world scenarios \citep{laidoui2023}. Game theory complements this approach by providing a framework to analyze both conflicts and cooperation, facilitating optimal outcomes for stakeholders. Together, these methods offer deeper insights and effective solutions to complex, multi-stakeholder challenges, which is the focus of this paper.

Most research on multi-objective charging scheduling, like many other optimization problems, relies on decision space search algorithms. For example, fuzzy genetic algorithms and fuzzy discrete particle swarm optimization have been used to coordinate the charging of PEVs \citep{hajforoosh2015}. Similarly, genetic algorithms, differential evolution, and particle swarm optimization have been applied to tackle charging scheme challenges \citep{wu2017}. In online EV charge scheduling, a heuristic algorithm combines a multi-commodity network flow model with a customized bisection search within a rolling horizon framework \citep{wu2020}. An improved multi-objective particle swarm optimization algorithm has also been developed to address EV charging scheduling \citep{yin2021}. However, these decision space search algorithms often lack optimality guarantees, which can be a significant limitation.

In contrast, few studies have explored multi-objective charging scheduling using criterion space search algorithms. For example, \cite{zakariazadeh2014} used the augmented $\epsilon$-constraint method, and \cite{das2020} proposed an augmented non-dominated $\epsilon$-constraint algorithm for multi-objective EV charging. These methods, however, are generally more suited to problems that prioritize a primary objective, as they focus on identifying a subset of non-dominated solutions rather than the entire set.

In addition to the $\epsilon$-constraint method, commonly used criterion space search methods include the weighted sum method, perpendicular search method, and augmented weighted Tchebycheff method. These methods operate by converting a multi-objective optimization problem into a series of single-objective problems, which are then solved sequentially. As a relatively novel criterion space search method, the balanced box method has demonstrated effectiveness. According to \cite{boland2015}, the balanced box method decomposes the criterion space and has been successfully applied to various bi-objective optimization problems \citep{wang2020, charkhgard2020, smet2023}. Building on this, some studies have proposed extensions of the balanced box method. For example, \cite{dai2018} introduced a two-stage approach that combines the balanced box method with the $\epsilon$-constraint method to solve a smaller set of single-objective integer linear programs (SOILPs).

While decision space search algorithms have been widely applied to multi-objective EV charging scheduling, they often lack guarantees of global optimality and robustness. Criterion space search methods, such as the $\epsilon$-constraint method, offer more precise identification of optimal solutions but are computationally intensive and tend to emphasize a single objective. In this paper, we extend the balanced box method and propose two enhanced strategies to solve bi-objective optimization problems more effectively. These strategies address the need for efficient, robust methods that provide high-quality solutions at lower computational costs while supporting collaborative scheduling among stakeholders in shared charging systems.

\subsection*{Summary}

EV charging scheduling is a complex problem that has attracted considerable research attention. Various models have been developed to optimize charging schedules in settings such as parking garages, fleet management systems, and shared charging environments. Conflicting objectives are often considered in EV charging management to balance competing priorities. To tackle these challenges, researchers have applied a range of optimization algorithms, including heuristic methods such as particle swarm optimization and genetic algorithms, and exact methods like the $\epsilon$-constraint method. While these methods have contributed to advances in the field, they often overlook robustness and fairness perspectives that are critical in real-world applications.

To address the above challenges, we propose a collaborative scheduling framework that integrates bi-objective optimization with game theory. This framework focuses on shared charging stations, balancing the interests of two companies sharing charging resources. In contrast to traditional heuristic-based approaches, we focus on criterion space search methods. Building on the Balanced Box Method, we developed the B3M to efficiently delineate the non-dominated frontier and provide a compact yet representative solution set. By combining B3M with cooperative bargaining, our approach effectively determines a final agreement point that ensures balanced benefits for both parties, thereby promoting fair and cooperative decision-making. This innovation addresses key limitations of existing methods and provides a robust, scalable solution that enhances collaboration in real-world shared charging station contexts.

\section{Bi-objective integer programming model} \label{Section2}
\subsection{Problem description}

\begin{figure*} [htbp]
	\centering
		\includegraphics[scale=.11]{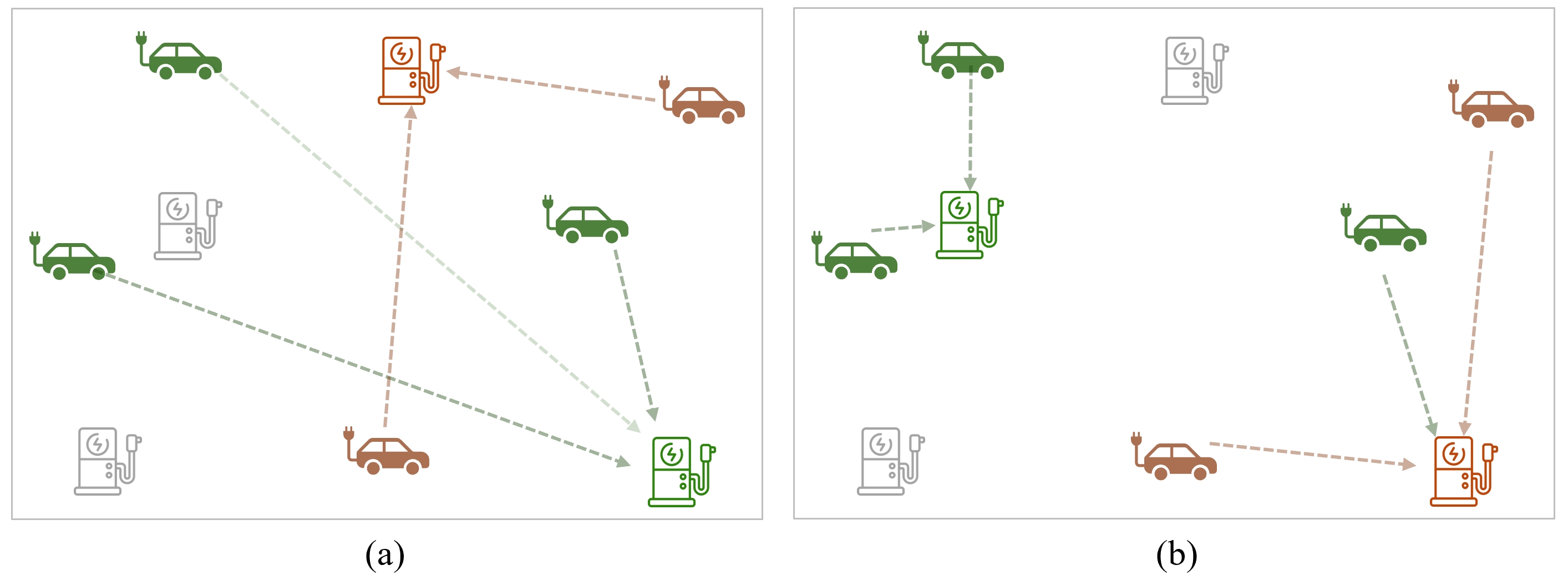}
    \caption{Illustration of collaborative charging scheduling problem}
	\label{FIG:illustration}
\end{figure*}

This study examines the collaborative scheduling problem where two fleet operators rent and share dedicated charging facilities from a third party to exclusively serve their EV fleets. The goal is to jointly optimize charging station rental decisions and charging schedules while minimizing operational costs and ensuring fair resource allocation, as illustrated in Fig. \ref{FIG:illustration}.

In this setup, the two companies (i.e., the fleet operators, represented in green and orange) have access to four charging stations available for rent. Fig. \ref{FIG:illustration} (a) depicts the non-collaborative scenario, where each charging station is exclusively used by a single company. In contrast, Fig. \ref{FIG:illustration} (b) illustrates the collaborative scenario, where EVs from both companies can utilize all charging stations within the coalition. By sharing charging resources, the collaboration is expected to lower costs, improve utilization, and provide greater scheduling flexibility.

The core challenge in this problem is jointly optimizing charging station rental decisions and charging schedules to minimize operational costs while ensuring efficient and fair resource allocation. Each company seeks to minimize its own costs, including station rental fees, electricity costs, energy consumption for traveling to the station, and waiting time due to scheduling constraints. These competing objectives naturally lead to a bi-objective optimization problem, where solutions form an efficient frontier balancing cost efficiency and fairness. To ensure equitable resource distribution, a cooperative bargaining mechanism is introduced, guiding companies toward a fair and cost-effective solution.

To guarantee feasible scheduling, the problem formulation incorporates multiple constraints: (i) charging continuity constraints, ensuring uninterrupted charging sessions, (ii) charger capacity constraints, limiting simultaneous usage per station, (iii) charging time constraints, defining valid scheduling windows, (iv) charging demand and capacity constraints, ensuring sufficient energy supply for each EV, and (v) charger rental constraints, ensuring each station is rented by only one company and used accordingly.

\subsection{Model setup and notations}

We consider a collaborative EV charging scenario in which multiple fleet operators coordinate the use of a set of third-party charging stations via a centralized scheduling platform. This structure reflects an emerging operational model in urban logistics and shared mobility systems, where vehicles from different companies access charging infrastructure under rental agreements and negotiated usage rules.

In such systems, the next day's operational plans and vehicle assignments are typically determined in advance. Fleet operators use telematics and battery monitoring systems to estimate charging needs and submit them to the platform prior to scheduling. Charging demand, including energy requirements, time windows, and location preferences, is assumed to be known, as it is typically derived from upstream planning processes such as route planning and battery monitoring. This reflects the scheduling-level focus of the proposed model. Based on these inputs, the platform determines which vehicles are granted access to which stations, during which time intervals, and under what pricing terms.

Charging station usage involves two types of cost: rental and electricity fees. A company that leases a station is responsible for the corresponding rental fee and electricity cost incurred when its own EVs are charged. Under centralized coordination, the use of a leased station may be shared: if a vehicle from another participating company is scheduled to use that station, it is charged a collaborative electricity price defined by the platform. Electricity prices may vary by user type, station, and time interval. To reflect physical capacity constraints, each charger is restricted to serving a single EV per time slot in the model.

The scheduling horizon is discretized into equal-length time intervals. Vehicles can only be assigned to contiguous time slots, and their total charging duration must be sufficient to meet their required energy. Moreover, we assume that once charging begins, it must be uninterrupted until the energy demand is satisfied. This reflects typical operational constraints, such as physical plug-locking mechanisms or platform-imposed charging continuity rules.

Given these assumptions and modeling setups in place, we introduce decision variables for charging station rental and EV scheduling. Specifically, $y_j^k$ indicates whether company $k$ rents charger $j$, while $x_{ij}^\tau$ represents whether EV $i$ is assigned to charger $j$ during time interval $\tau$. In addition, time-related variables such as $t_{i}^{\rm s}$ and $t_{i}^{\rm f}$ are introduced to track the start and end times of each EV’s charging session, ensuring proper scheduling coordination. The notations used in this study, including decision variables, auxiliary variables, parameters, and sets, are summarized in Table \ref{tbnotation} for reference.
To formally define the optimization problem, we distinguish between predefined input parameters and optimization decision variables.
The input parameters include the set of EVs $I$, chargers $J$, companies $K$, and time intervals $T$, along with parameters such as charging rates, energy fees, rental cost, and vehicle-specific charging demands. The decision variables include charging station rental decision ($y_j^k$) and EV charging schedules ($x_{ij}^\tau$), while auxiliary variables facilitate scheduling coordination and maintain model tractability. These elements collectively define the bi-objective integer programming model presented in Section \ref{modelformulation}.

\begin{table}[width=.9\linewidth,cols=2,pos=h]
\caption{Mathematical notation}\label{tbnotation}
\begin{tabular*}{\tblwidth}{@{} LL@{} }
\toprule
\textbf{Sets:} &                    \\
$I$ & Set of EVs, $i \in I$\\
$J$ & Set of chargers, $j \in J$\\
$K$ & Set of companies, $k \in K$\\
$T$ & Set of time intervals, $\tau \in T$\\
$I_k$ & Set of EVs of company $k$, $I_k \subset I$\\
\textbf{Parameters:} &      \\
${\rm f}_j^k$ & Fixed rental fee of charger $j$ for company $k$\\
${\rm c}_{\rm e,c}^{j,\tau}$ & Unit energy fee for collaborative EVs charging at charger $j$ in time interval $\tau$\\
${\rm c}_{\rm e,o}^{j,\tau}$ & Unit energy fee for own EVs charging at charger $j$ in time interval $\tau$\\
${\lambda}_i$ & Value of time of EV $i$\\
$\varepsilon_{ij}$ & Charging rate of EV $i$ at charger $j$\\
${w}_{ij}$ & Energy cost for EV $i$ travels to charger $j$\\
$\left[e_i, l_i\right]$ & Available charging time window of EV $i$\\
$\left[L_i, H_i\right]$ & Maximum and minimal charging amount of EV $i$\\
\textbf{Decision variables:} &   \\
$x_{ij}^\tau$ & 1 if charger $j$ serves EV $i$ in time interval $\tau$, 0 otherwise\\
$y_j^k$ & 1 if charger $j$ is rented by company $k$, 0 otherwise\\
\textbf{Auxiliary variables:} &   \\
$u_{ij}^{\pi,k}$ & Binary variable introduced to linearize the product of $x_{ij}^\tau$ and $y_j^k$.\\
$t_{i}^{\rm s}$, $t_{i}^{\rm f}$ & Start and end charging time of EV $i$\\
$\underline{x_{ij}^{\tau}}$, $\overline{x_{ij}^{\tau}}$ & 1 if EV $i$ starts (concludes) charging at charger $j$ during time interval $\tau$, 0 otherwise.\\
\bottomrule
\end{tabular*}
\end{table}

With these definitions in place, the following section presents the mathematical formulation of the problem, ensuring that all rental and scheduling decisions satisfy the defined constraints while balancing cost efficiency and fairness in a bi-objective framework.

\subsection{Model formulation} \label{modelformulation}

The problem is formulated as a bi-objective integer programming model, where all decision variables are integers. This formulation jointly optimizes station rental and scheduling while balancing cost efficiency and fairness. The model comprises two distinct objective functions, each representing the cost-minimization goal of a company.

\vspace{2 mm}
\textbf{\textit{(I) Objective function}}

The objective function is formulated as:

\vspace{-\topsep}
\begin{equation} \label{Obj}
\min_{(x_{ij}^\tau, y_j^k)\in \mathscr{X}} \left( z_k(x_{ij}^\tau, y_j^k) \mid k = 1, 2 \right).
\end{equation}

In Eq. (\ref{Obj}), $z_k (x_{ij}^\tau, y_j^k)$ represents the costs incurred by company $k$. These costs are comprised of several sources, including rental fees, electricity charges (from EV charging), energy consumption (from travel), and time costs. It is defined as follows:

\vspace{-\topsep}
\begin{equation} \label{Objk1}
\begin{split}
	z_k (x_{ij}^\tau, y_j^k) = &\sum_{j \in J}{\rm f}_j^k y_j^k + \sum_{i \in I_k}\sum_{j \in J}\sum_{\tau \in T} \left(\sum_{{k'} \in K,{k'}\neq k}{\rm c}_{\rm e,c}^{j,\tau}y_j^{k'} + {\rm c}_{\rm e,o}^{j,\tau} y_j^{k}\right)\varepsilon_{ij} x_{ij}^\tau \\
    & + \sum_{i\in I_k}\sum_{j \in J}\sum_{\tau \in T} {w}_{ij}\underline{x_{ij}^\tau} + \sum_{i \in I_k}{\lambda}_i \left(t_i^{\rm s}-e_i\right),
\end{split}
\end{equation}

\noindent where the first term represents the rental fee, which is determined by the fixed unit cost, represented by ${\rm f}_j^k$, for the rented charging stations. The second term denotes the charging fee. The unit charging fee for collaborative EVs and own EVs is represented by ${\rm c}_{\rm e,c}^{j,\tau}$ and ${\rm c}_{\rm e,o}^{j,\tau}$, respectively. Moreover, the variable $\varepsilon_{ij}$ represents the charging rate. The third term is the energy consumption cost from the EV to the charging station, with the energy consumption matrix represented by ${w}_{ij}$. The final term represents the opportunity cost associated with waiting time for charging. This can be calculated by multiplying the drivers' value of time (${\lambda}_i$) by the waiting time ($t_i^{\rm s}-e_i$). It should be noted that electricity fees may vary depending on location and time, and may differ between charging the company's own EVs and those of the other company.

In this model, the two objectives are addressed separately, rather than being combined into a single weighted sum. This approach yields a \textbf{non-dominated frontier} also referred to as the \textbf{efficient frontier}, wherein the advancement of one objective entails a detriment to the other. The efficient frontier offers a means of reconciling the two companies' goals, a process that will be further elucidated in Section \ref{Section3}.

\vspace{2 mm}
\textbf{\textit{(II) Charging continuity constraints}}

To ensure uninterrupted charging, two additional auxiliary variables, denoted by $\underline{x_{ij}^{\tau}}$ and $\overline{x_{ij}^{\tau}}$, are introduced. If the value of $\underline{x_{ij}^{\tau}}$ is equal to 1, this signifies that EV $i$ starts the charging process at charger $j$ within the specified time interval $\tau$, as defined in: 

\vspace{-\topsep}
\begin{equation} \label{xstart}
	\underline{x_{ij}^{\tau}} \geq 
 \begin{cases}
 x_{ij}^{\tau}-x_{ij}^{\tau-1}, &\tau=2,3,4\dots, \mathscr{T}\\
 x_{ij}^{\tau}, &\tau=1.
 \end{cases}
 \forall i\in I, j \in J
\end{equation}

Similarly, if the value of $\overline{x_{ij}^{\tau}}$ is equal to 1, this denotes that EV $i$ has completed the charging process at charger $j$ within the time interval $\tau$, as given by:

\vspace{-\topsep}
\begin{equation} \label{xend}
	\overline{x_{ij}^{\tau}} \geq 
 \begin{cases}
 x_{ij}^{\tau}-x_{ij}^{\tau+1}, &\tau=1,2,3\dots, \mathscr{T}-1\\
 x_{ij}^{\tau}, &\tau=\mathscr{T}.
 \end{cases}
 \forall i\in I, j \in J
\end{equation}

\vspace{2 mm}
\textbf{\textit{(III) Charging location constraints}}

Each EV is constrained to charging at a single charger only once, and these limitations can be expressed as:

\vspace{-\topsep}
\begin{equation} \label{xstart1}
	\sum_{\tau \in T}\sum_{j \in J}\underline{x_{ij}^{\tau}} =1, \forall i \in I,
\end{equation}

\vspace{-\topsep}
\begin{equation} \label{xend1}
	\sum_{\tau \in T}\sum_{j \in J}\overline{x_{ij}^{\tau}} =1, \forall i \in I,
\end{equation}

\noindent and

\vspace{-\topsep}
\begin{equation} \label{xstartend}
	\sum_{\tau \in T}\underline{x_{ij}^{\tau}}=\sum_{\tau \in T}\overline{x_{ij}^{\tau}}, \forall i \in I, j \in J.
\end{equation}

\vspace{2 mm}
\textbf{\textit{(IV) Charger capacity constraints}}

Each charger is limited to charging one EV per time interval, as expressed by:

\vspace{-\topsep}
\begin{equation} \label{x1}
	\sum_{i \in I}x_{ij}^{\tau}\leq 1, \forall j \in J, \tau \in T.
\end{equation}

\vspace{2 mm}
\textbf{\textit{(V) Charging time constraints}}

The start charging time $t_{i}^{\rm s}$ and end charging time $t_{i}^{\rm f}$ are also auxiliary variables that can be represented by:

\vspace{-\topsep}
\begin{equation} \label{tstart}
	t_{i}^{\rm s} = \sum_{\tau \in T}\sum_{j \in J} \underline{x_{ij}^{\tau}} (\tau-1), \forall i \in I,
\end{equation}

\vspace{-\topsep}
\begin{equation} \label{tend}
	t_{i}^{\rm f} =\sum_{\tau \in T}\sum_{j \in J} \overline{x_{ij}^{\tau}} \tau, \forall i \in I,
\end{equation}

\noindent and

\vspace{-\topsep}
\begin{equation} \label{xend2}
	\sum_{\tau \in T}\sum_{j \in J}x_{ij}^{\tau} = t_{i}^{\rm f} - t_{i}^{\rm s},\forall i\in I.
\end{equation}

To clarify the representation of the charging duration, the start charging time $t_{i}^{\rm s}$ is defined to occur between 0 and $\mathscr{T}-1$, while the end charging time $t_{i}^{\rm f}$ occurs between 1 and $\mathscr{T}$. 
The total charging duration is calculated as the difference between the end charging time $t_i^{\rm f}$ and the start charging time $t_i^{\rm s}$. For example, if charging takes place during the interval $\tau$, the start time is $\tau-1$, and the end time is $\tau$, resulting in a charging duration of one time unit.

The start and end of charging are constrained by the predetermined time windows for each EV:

\vspace{-\topsep}
\begin{equation} \label{tstartendTW}
	e_i \leq t_{i}^{\rm s} \leq t_{i}^{\rm f} \leq l_i, \forall i\in I.
\end{equation}

\vspace{2 mm}
\textbf{\textit{(VI) Charging demand and capacity constraints}}

The minimum and maximum requirements limit the amount of charge for each EV:

\vspace{-\topsep}
\begin{equation} \label{Demandi}
L_i\leq \sum_{j \in J}\sum_{\tau \in T}\varepsilon_{ij} x_{ij}^\tau \leq H_i, \forall i\in I, 
\end{equation}

\noindent where the charging rate, denoted as $\varepsilon_{ij}$, corresponds to the charging speed of charger $j$ for EV $i$. It is important to note that this charging rate is influenced not only by the characteristics of both chargers and EVs. $L_i$ and $H_i$ are the lower and upper bounds of the charging demand of EV $i$.

\vspace{2 mm}
\textbf{\textit{(VII) Charger rental constraints}}

Each charger can be rented by a maximum of one company, as constrained by:

\vspace{-\topsep}
\begin{equation} \label{zjk}
	\sum_{k \in K}y_j^k \leq 1, \forall j \in J.
\end{equation}

Furthermore, EVs can only be served by a charger if it is rented, expressed as:

\vspace{-\topsep}
\begin{equation} \label{xijzjk}
	x_{ij}^\tau \leq \sum_{k \in K}y_j^k, \forall i \in I, j\in J, \tau \in T.
\end{equation}

With these constraints, the problem formulation is complete. 

\vspace{2 mm}

To simplify the solution process for the proposed model, the objective function is linearized. An auxiliary binary variable, $u_{ij}^{\tau,k}$, is introduced to replace the nonlinear term $x_{ij}^\tau y_j^{k}$. This transformation enables the use of standard integer linear programming techniques. The following constraints are introduced to ensure the linearization accurately captures the relationship between $x_{ij}^\tau$, $y_j^k$, and $u_{ij}^{\tau,k}$:

\begin{equation}
\begin{aligned}
    u_{ij}^{\tau,k} &\leq x_{ij}^\tau, \\
    u_{ij}^{\tau,k} &\leq y_j^k, \\
    u_{ij}^{\tau,k} &\geq x_{ij}^\tau + y_j^k - 1.
\end{aligned}
\end{equation}

Consequently, the objective function can be expressed as follows:

\vspace{-\topsep}
\begin{equation} \label{Objk2}
 \begin{split}
	z_k (x_{ij}^\tau, y_j^k)  = & \sum_{j \in J}{\rm f}_j^k y_j^k + \left(\sum_{i \in I}\sum_{j \in J}\sum_{\tau \in T} \sum_{{k'} \in K,{k'}\neq k}{\rm c}_{\rm e,c}^{j,\tau} \varepsilon_{ij} u_{ij}^{\tau,k'} + \sum_{i \in I}\sum_{j \in J}\sum_{\tau \in T} {\rm c}_{\rm e,o}^{j,\tau} \varepsilon_{ij} u_{ij}^{\tau,k}\right) \\
 & + \sum_{i\in I_k}\sum_{j \in J}\sum_{\tau \in T} {w}_{ij}\underline{x_{ij}^\tau} + \sum_{i \in I_k}{\lambda}_i \left(t_i^{\rm s}-e_i\right).
 \end{split}
\end{equation}

To summarize, the collaborative scheduling problem is modeled by Eqs. (1) and (3)–(17), which is inherently a \textbf{bi-objective integer linear programming model}. The model incorporates binary decision variables $x_{ij}^\tau$ and $y_j^k$, as well as auxiliary binary variables $\underline{x_{ij}^\tau}$, $\overline{x_{ij}^\tau}$ and $u_{ij}^{\tau,k}$ alongside auxiliary integer variables $t_i^{\rm s}$ and $t_i^{\rm f}$. By minimizing the costs for each company, the two objective functions address the inherent conflicts in collaborative scheduling. The model also incorporates key constraints, including charging continuity, location, time, demand, and capacity constraints, as well as charger capacity and rental constraints, ensuring the practicality and feasibility of the scheduling framework.

\section{Solution methods} \label{Section3}

This section presents the methods used to solve the above problem, which is divided into two parts. The first part uses criterion space search methods, specifically the development of the B3M, to identify the efficient frontier for the bi-objective function. The second part applies cooperative bargaining techniques to select the optimal solution from the solutions generated in the first part. The B3M represents a significant enhancement of the balanced box method, addressing its limitations in handling computational efficiency and solution representativeness. Unlike the original balanced box method, which focuses solely on accurately delineating the exact efficient frontier, the B3M introduces a novel mechanism to generate an efficient frontier. This mechanism strategically filters out closely positioned solutions, thereby reducing computational complexity while preserving the diversity and representativeness of the solution set.  Once the efficient frontier has been established, cooperative bargaining methods are implemented to select the final optimal solution, ensuring a balanced and collaborative outcome. By integrating these innovations, the B3M not only maintains the strengths of the original balanced box method but also significantly enhances its applicability to complex bi-objective optimization problems.

Before introducing the solution approaches, it is essential to delineate the fundamental concepts. The notation $(x_{ij}^{\tau}, y_j^k)$ is used to represent a combination of decision variables, where $x_{ij}^{\tau}$ and $y_j^k$ are individual decision variables. The objective function values associated with this solution are represented by $z_1(x_{ij}^{\tau}, y_j^k)$ and $z_2(x_{ij}^{\tau}, y_j^k)$, which correspond to the first and second objective functions, respectively. For simplicity, we use the shorthand notations $z_1(\mathscr{x})$ and $z_2(\mathscr{x})$ throughout this study, where $\mathscr{x}$ represents the combination of decision variables $(x_{ij}^{\tau}, y_j^k)$.

For a given solution $\mathscr{x}$, the corresponding objective function value $z$ is:

\vspace{-\topsep}
\begin{equation} \label{Objvalue}
z(\mathscr{x}) = \left( z_1(\mathscr{x}), z_2(\mathscr{x}) \right).
\end{equation}

The set of objective function values corresponding to the solution set is typically referred to as the objective space:

\vspace{-\topsep}
\begin{equation} \label{Objset}
\mathcal{Z} = \left\{ z(\mathscr{x}) \mid \mathscr{x} \in \mathscr{X} \right\},
\end{equation}

\noindent where $\mathscr{X}$ is the set of solutions, and the objective space is represented by $\mathcal{Z}$. The latter contains the values of all solutions with respect to the objective function.

The notations used for all the solution approaches can be found in Table \ref{tbnotation1}.

\begin{table}[width=.9\linewidth,cols=2,pos=h]
\caption{Notations for solution approaches}\label{tbnotation1}
\begin{tabular*}{\tblwidth}{@{} p{2cm} p{12cm} @{}}
\toprule
Symbol & Description \\
\midrule
$\mathscr{X}, \mathcal{Z}$ & Set of solutions and corresponding objective function values, where $\mathscr{x} \in \mathscr{X}$ and $z \in \mathcal{Z}$ \\
$\tilde{\mathscr{X}}, \tilde{\mathcal{Z}}$ & Reduced set of solutions and corresponding objective function values, where $\tilde{\mathscr{X}} \subset \mathscr{X}$ and $\tilde{\mathcal{Z}} \subset \mathcal{Z}$ \\
$\mathscr{x}, z$ & A specific solution and its corresponding objective function value, where $z=z(\mathscr{x})$, and $\mathscr{x}$ represents the combination of decision variables $(x_{ij}^{\tau}, y_j^k)$ \\
$\mathscr{x}^*, z^*$ & Final solution and its corresponding objective function value \\
$z^T, z^B$ & The top and bottom points of the efficient frontier \\
$z_k$ & Objective function $k$ \\
$z_1, z_2$ & The first and second objective functions, where $z(\mathscr{x}) = \left( z_1(\mathscr{x}), z_2(\mathscr{x}) \right)$ \\
$\mathcal{C}_k^{\rm min}$ & Cost threshold of company $k$ \\
$z^{\rm Non}$ & Non-collaborative point\\
$z^t, z^b$ & Points to define the top and bottom rectangles \\
$z^{\rm n,1}, z^{\rm n,2}$ & New solutions found in the bottom rectangle and top rectangle\\
$R(z^T,z^B)$ & The rectangle defined by the top-left point $z^T$ and the bottom-right point $z^B$ \\
$R^T, R^B$ & Top and bottom rectangles\\
$\mathcal{Z}^{\rm e}$ & The set of objective function values corresponding to the existing solutions, where $z^{\rm e} \in \mathcal{Z}^{\rm e}$ \\
$z^{\rm e}$ & Existing solution\\
$z^{\rm n}$ & New solution\\
$\sigma_1, \sigma_2$ & Tolerance margins for objectives 1 and 2 \\
$\epsilon$ & Tolerance range, expressed as a percentage \\
$\zeta$ & A small constant \\
\bottomrule
\end{tabular*}
\end{table}

\subsection{Balanced bounding box methods} \label{Section3.1}

The balanced bounding box method (B3Ms) is inspired by the \textbf{balanced box method}, which has been recognized for its ability to efficiently identify the exact efficient frontier compared to other criterion space exploration techniques \citep{boland2015}. The balanced box method systematically explores rectangles within the criterion space, with the flowchart and main steps provided in the Appendix. A detailed description of its procedure is available in \citet{boland2015}.

\vspace{2 mm}

\noindent \textbf{\textit{Definition 1: Participation constraints in collaborative optimization}}

\textit{
A participation constraint ensures that a company’s cost in a collaborative solution does not exceed its cost in a non-collaborative scenario. In this study, we introduce a threshold condition guaranteeing that every company benefits from collaboration compared to its independent operation. This constraint does not define the entire feasible solution space but serves as a necessary condition to exclude unreasonably unfavorable outcomes and ensure voluntary participation in the collaboration.}

\subsubsection{Incorporating participation constraints for collaboration}

To ensure efficiency in collaborative optimization problems, the \textbf{B3Ms} framework imposes participation constraints that restrict solutions where a company's cost would exceed its non-collaborative scenario. Specifically, we introduce a threshold constraint \citep{fernandez2016,zhou2024}, ensuring that each company’s cost remains below a predefined benchmark derived from its independent operation. This benchmark, denoted as $z^{\rm Non}$, serves as an upper bound for the objectives in collaboration and is formulated as:

\vspace{-\topsep}
\begin{equation} \label{NonPk}
	z_k \leq z_k^{\rm Non}, \forall k \in K.
\end{equation}

\begin{figure*} [htbp]
	\centering
		\includegraphics[scale=.6]{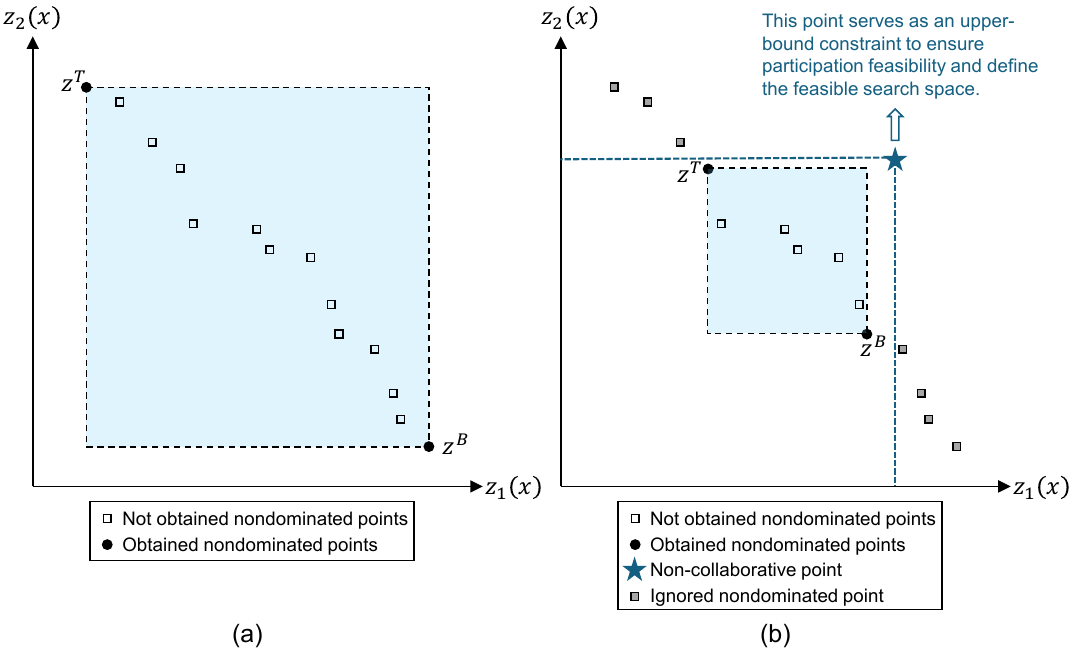}
    \caption{Balanced box method: (a) unbounded search space, (b) incorporating non-collaborative point}
	\label{FIG:0}
\end{figure*}

This condition effectively constrains the search space, ensuring that only solutions aligning with collaborative goals are considered. As illustrated in Fig. \ref{FIG:0}, the initial search space is constrained by the participation constraints, which use the non-collaborative point as a reference. This ensures that only solutions where both companies benefit from collaboration are considered feasible. Rather than directly influencing the determination of the efficient frontier, these constraints delineate the feasible region for collaboration, excluding solutions that result in worse outcomes than the non-collaborative scenario to guarantee voluntary participation.

Given these constraints, the initial and final points of the efficient frontier ($z^T$ and $z^B$) are determined by solving the following two equations:

\vspace{-\topsep}
\begin{equation} \label{zTnon}
	z^T:=\mathop{\rm{lexmin}}_{\mathscr{x}\in \mathscr{X}} \left\{z_1(\mathscr{x}),z_2(\mathscr{x}):z(\mathscr{x})\in R((0,z^{\rm Non}_1),(0,z^{\rm Non}_2))\right\}
\end{equation}

\noindent and

\vspace{-\topsep}
\begin{equation} \label{zBnon}
	z^B:=\mathop{\rm{lexmin}}_{\mathscr{x}\in \mathscr{X}} \left\{z_2(\mathscr{x}),z_1(\mathscr{x}):z(\mathscr{x})\in R((0,z^{\rm Non}_1),(0,z^{\rm Non}_2))\right\},
\end{equation}

\noindent respectively. In Eqs. \eqref{zTnon} and \eqref{zBnon}, "lexmin" denotes a lexicographic optimization approach, wherein the objectives are minimized in a sequential manner. This entails that each objective must be minimized individually. Eq. \eqref{zTnon} initially minimizes $z_1(\mathscr{x})$, subsequently minimizing $z_2(\mathscr{x})$ while treating the objective values of $z_1(\mathscr{x})$ as constraints. Similarly, Eq. \eqref{zBnon} applies the fundamental principle. Note that each point requires two processing stages to obtain the complete solution.

It should be noted that in defining the initial search space, we incorporate participation constraints using the non-collaborative point as a reference rather than considering an unbounded space. This ensures that only solutions where all companies benefit from collaboration are retained, effectively filtering out infeasible solutions and significantly reducing computational time.

Based on the identified initial and final points of the efficient frontier, the initial rectangle is defined as $R(z^T,z^B)$. The process of exploring rectangles involves dividing a rectangle into two smaller regions: the upper and lower rectangles ($R^T$ and $R^B$). These two rectangles are defined by the points $z^T$, $z^t$, $z^b$, and $z^B$, where $z^t=(z^B_1, (z^T_2+z^B_2)/2)$ and $z^b=(z^T_1, (z^T_2+z^B_2)/2)$. The original rectangle is divided horizontally along the $z_2(\mathscr{x})$ axis and then removed from the rectangle set. New non-dominated points are identified by solving lexicographic optimization problems within these regions, and new rectangles are added to update the rectangle set. Smaller rectangles created during the process are handled in the same way. This iterative process continues until the rectangle set is empty, at which point all non-dominated points are outputted. Further details of the two B3Ms will be introduced later.

\subsubsection{The fundamental concept of B3Ms} \label{section3.1.2}

The fundamental concept underlying the B3Ms is the exclusion of solutions that are close to one another. This ensures that the frontier generation process works efficiently and also speeds up computational time. The B3Ms optimize the solution framework to identify partial non-dominated solutions that are capable of representing the entire solution set. By excluding close solutions, the method avoids potential challenges in decision-making, such as difficulties in distinguishing among nearly equivalent solutions. Furthermore, it robustly discovers and represents the diversity that exists within the solution space, ensuring stable performance even when confronted with complex structures.

The reduced solution set obtained via the B3Ms is represented by $\tilde{\mathscr{X}}$. This is a subset of the solution set, denoted as $\tilde{\mathscr{X}} \subset \mathscr{X}$. In this context, the notation $\tilde{\mathcal{Z}}$ denotes the reduced set of objective function values derived from $\tilde{\mathscr{X}}$, defined as follows:

\vspace{-\topsep}
\begin{equation} \label{Objreducedset}
\tilde{\mathcal{Z}} = \left\{ z(\mathscr{x}) \mid \mathscr{x} \in \tilde{\mathscr{X}} \right\},
\end{equation}

\noindent where a reduced set of objective function values is represented by $\tilde{\mathcal{Z}} \subset \mathcal{Z}$, corresponding to the reduced solution set $\tilde{\mathscr{X}}$.

\vspace{2 mm}
\noindent \textbf{\textit{Definition 2: Relative closeness}}

\textit{We define a new solution $z^{\rm n}$ as being very close if it closely resembles an existing solution $z^{\rm e}$, where $z^{\rm e}$ belongs to the set $\mathcal{Z}^{\rm e}$. The degree of closeness is characterized by two distinct scenarios:}

\begin{itemize}
    \setlength{\itemsep}{0pt}
    \setlength{\parsep}{0pt}
    \setlength{\parskip}{0pt}
    \item \textbf{\textit{Definition 2.1: Strict closeness}} \textit{$z^{\rm n}\approx \mathcal{Z}^{\rm e}$ is defined if for all $z^{\rm e}\in \mathcal{Z}^{\rm e}$, the conditions $z^{\rm e}_1 - \sigma_1 \leq z^{\rm n}_1\leq z^{\rm e}_1 + \sigma_1$ and $z^{\rm e}_2 - \sigma_2 \leq z^{\rm n}_2 \leq z^{\rm e}_2+\sigma_2$ are simultaneously satisfied.}
    \item \textbf{\textit{Definition 2.2: Relaxed closeness}} \textit{$z^{\rm n}\sim \mathcal{Z}^{\rm e}$ is defined if for all $z^{\rm e}\in \mathcal{Z}^{\rm e}$, either the condition $z^{\rm e}_1 - \sigma_1 \leq z^{\rm n}_1\leq z^{\rm e}_1 + \sigma_1$ or $z^{\rm e}_2 - \sigma_2 \leq z^{\rm n}_2 \leq z^{\rm e}_2+\sigma_2$ is met.}
\end{itemize}

\noindent where $\sigma_1$ and $\sigma_2$ are tolerance margins for objectives 1 and 2, indicating the allowable range for closeness.

\begin{figure*} [htbp]
	\centering
		\includegraphics[scale=.2]{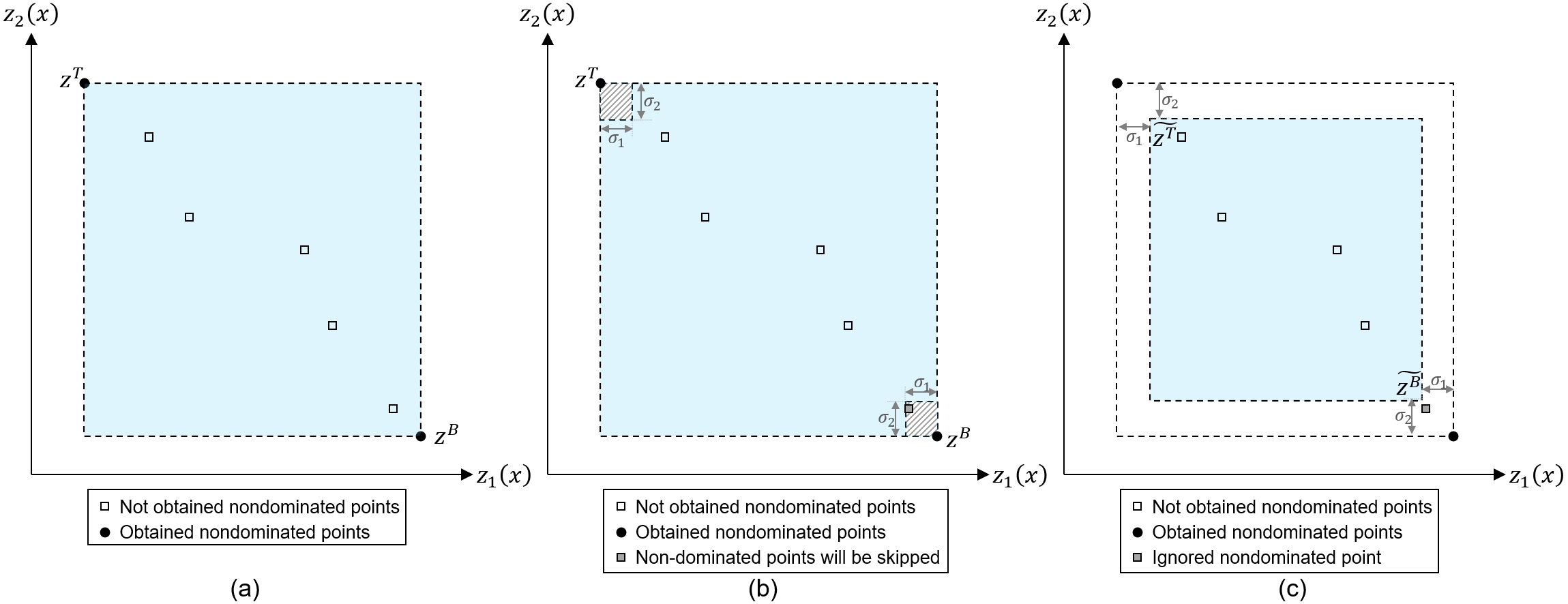}
    \caption{Balanced box methods (a) original, (b) B3M with $z^{\rm n}\approx \mathcal{Z}^{\rm e}$, (c) B3M with $z^{\rm n}\sim \mathcal{Z}^{\rm e}$}
	\label{FIG:1}
\end{figure*}

Two B3Ms were developed, each grounded in the concept of a respective closeness. The first, designated B3M1, is defined as $z^{\rm n}\approx \mathcal{Z}^{\rm e}$, while the second, designated B3M2, is defined as $z^{\rm n}\sim \mathcal{Z}^{\rm e}$. Fig. \ref{FIG:1} illustrates all the balanced box methods, including the original balanced box method, B3M1, and B3M2. The search area of interest is shown in blue. This is after the initial and final points of the efficient frontier have been determined.

As shown in Fig. \ref{FIG:1}, the search area of the B3M1 excludes only the small rectangles in the upper-left and lower-right corners. The B3M2, on the other hand, restricts the search area to the inner rectangle only.

\subsubsection{Balanced bounding box method 1}

\textit{The Balanced Bounding Box Method 1 (B3M1) is based on the strict closeness condition defined in Definition 2.1, where a new solution, denoted $z^{\rm n}$, is deemed to be sufficiently close to the existing solution set, denoted $\mathcal{Z}^{\rm e}$, to determine whether it is sufficiently close. This approach ignores two minor rectangular regions (the upper-left and lower-right corners) to streamline the search space. Once a non-dominated point has been identified, the method evaluates its closeness according to the strict conditions of Definition 2.1 and terminates the search in the current region if the conditions are met.}

\begin{figure*} [htbp]
	\centering
		\includegraphics[scale=.6]{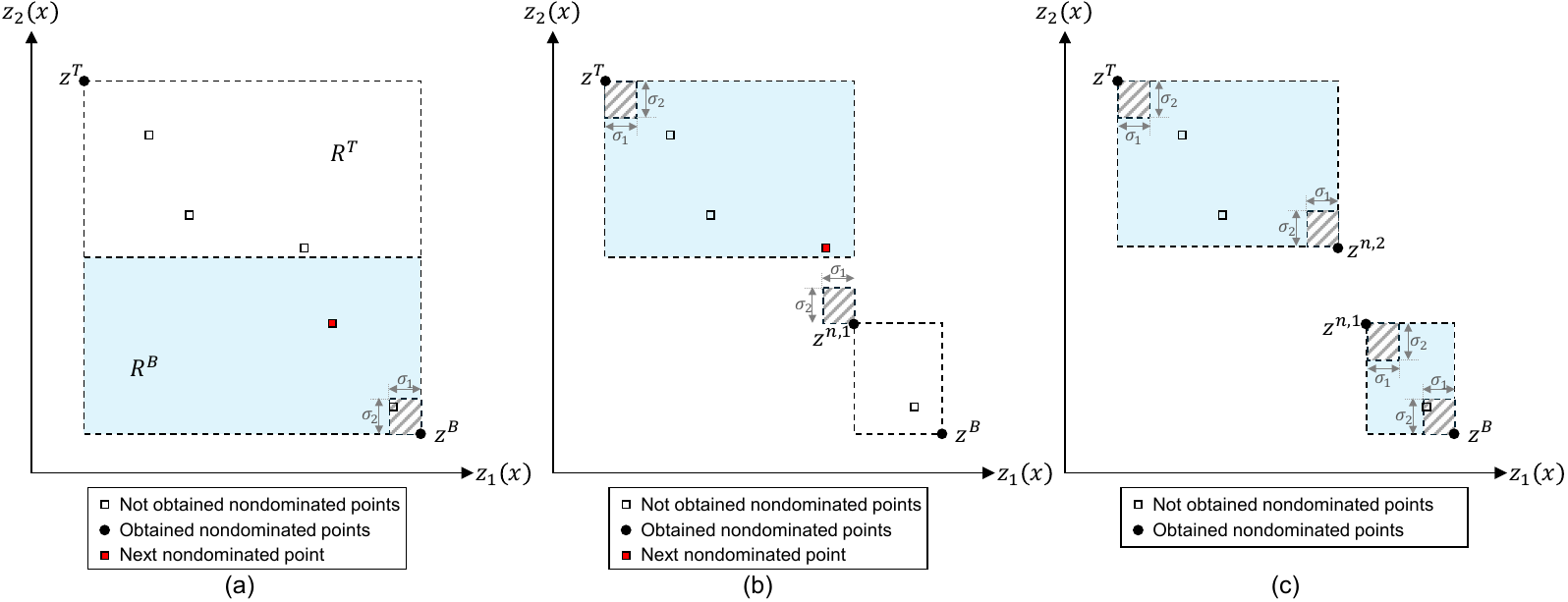}
    \caption{B3M1, $z^{\rm n}\approx \mathcal{Z}^{\rm e}$, (a) investigating $R^B$, (b) investigating $R^T$, (c) new rectangles}
	\label{FIG:M1_1}
\end{figure*}

\begin{figure*} [htbp]
	\centering
		\includegraphics[scale=.6]{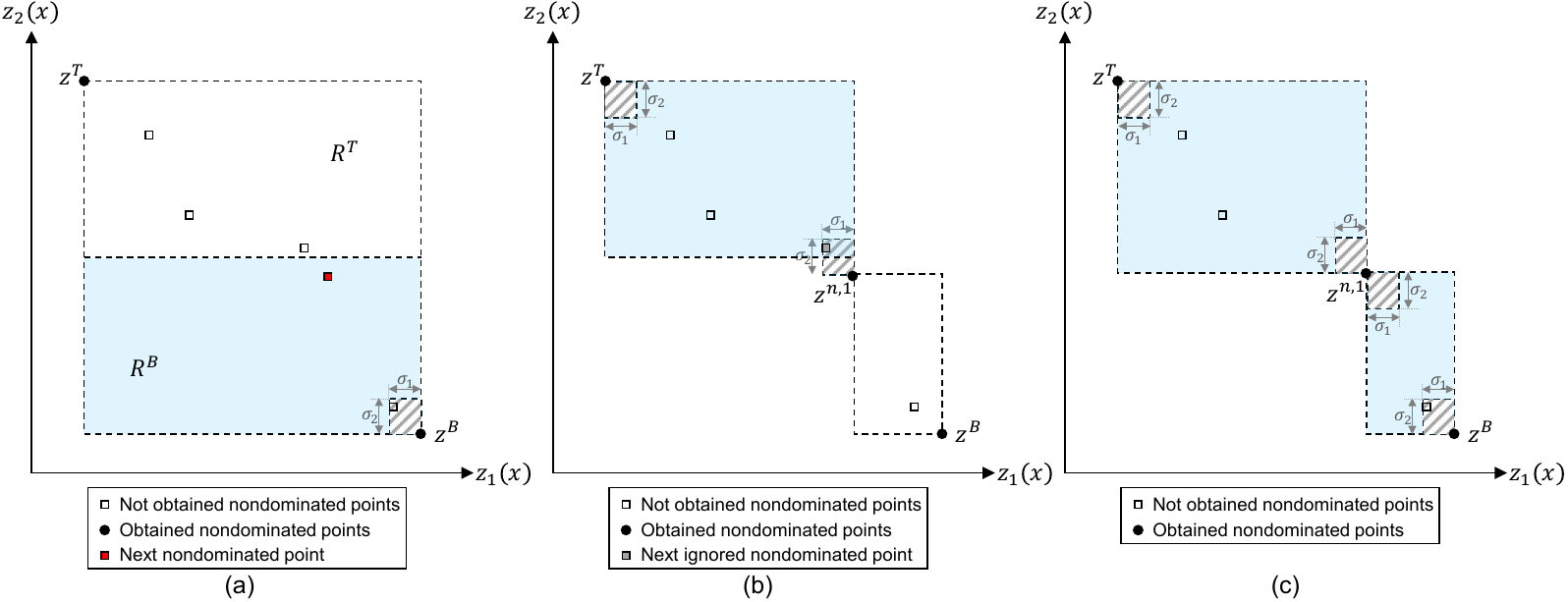}
    \caption{B3M1, $z^{\rm n}\approx \mathcal{Z}^{\rm e}$, \textbf{special scenario},  (a) investigating $R^B$, (b) investigating $R^T$, (c) new rectangles}
	\label{FIG:M1_2}
\end{figure*}

As shown in Fig. \ref{FIG:1} (b), for each rectangle, only the small rectangles in the upper-left and lower-right corners, marked by hatched areas, are ignored. The detailed steps of the method, including the comparison process in the lower and upper rectangles, are illustrated in Figs. \ref{FIG:M1_1} and \ref{FIG:M1_2}.

Fig. \ref{FIG:M1_1} (a) illustrates the search in the lower rectangle $R^B$, where the new solution is compared to the lower-right point. If the new solution lies outside the hatched area, it is recorded as a new non-dominated point; otherwise, no new rectangles are identified. Fig. \ref{FIG:M1_1} (b) shows the search in the upper rectangle $R^T$, where the new solution is compared with the upper-left point and the newly obtained solution $z^{\rm n,1}$ from the lower rectangle, or the lower-right point if $z^{\rm n,1}$ does not exist. As illustrated in Fig. \ref{FIG:M1_1} (b), since the new solution is not close to either $z^T$ or $z^{\rm n,1}$, it is recorded as a new non-dominated point. After searching the lower and upper rectangles, two new non-dominated points are obtained, similar to the original balanced box method. Consequently, two new rectangles are identified, as shown in Fig. \ref{FIG:M1_1} (c).

Nevertheless, there is a special situation where the solution found in the upper rectangle is close to the solution found in the lower rectangle, as illustrated in Fig. \ref{FIG:M1_2} (b). The new solution is close to $z^{\rm n,1}$ and should therefore be ignored. The non-dominated points obtained in this case are $z_T$, $z_B$, and $z^{\rm n,1}$, resulting in the two new rectangles depicted in Fig. \ref{FIG:M1_2} (c). Notably, the search area in Fig. \ref{FIG:M1_2} (c) is larger compared to the scenario where two new non-dominated points are identified, as shown in Fig. \ref{FIG:M1_1} (c). In this situation, the rectangle $R(z^T, \overline{z}^2)$ and the rectangle $R(\overline{z}^1, z^B)$ in Fig. \ref{FIG:M1_1} (c) can also represent the new search area of the original balanced box method. Consequently, this special case results in a larger search area due to the exclusion of some solutions, which we speculate may lead to an increase in computational time.

\begin{figure*} [htbp]
	\centering
		\includegraphics[scale=0.58]{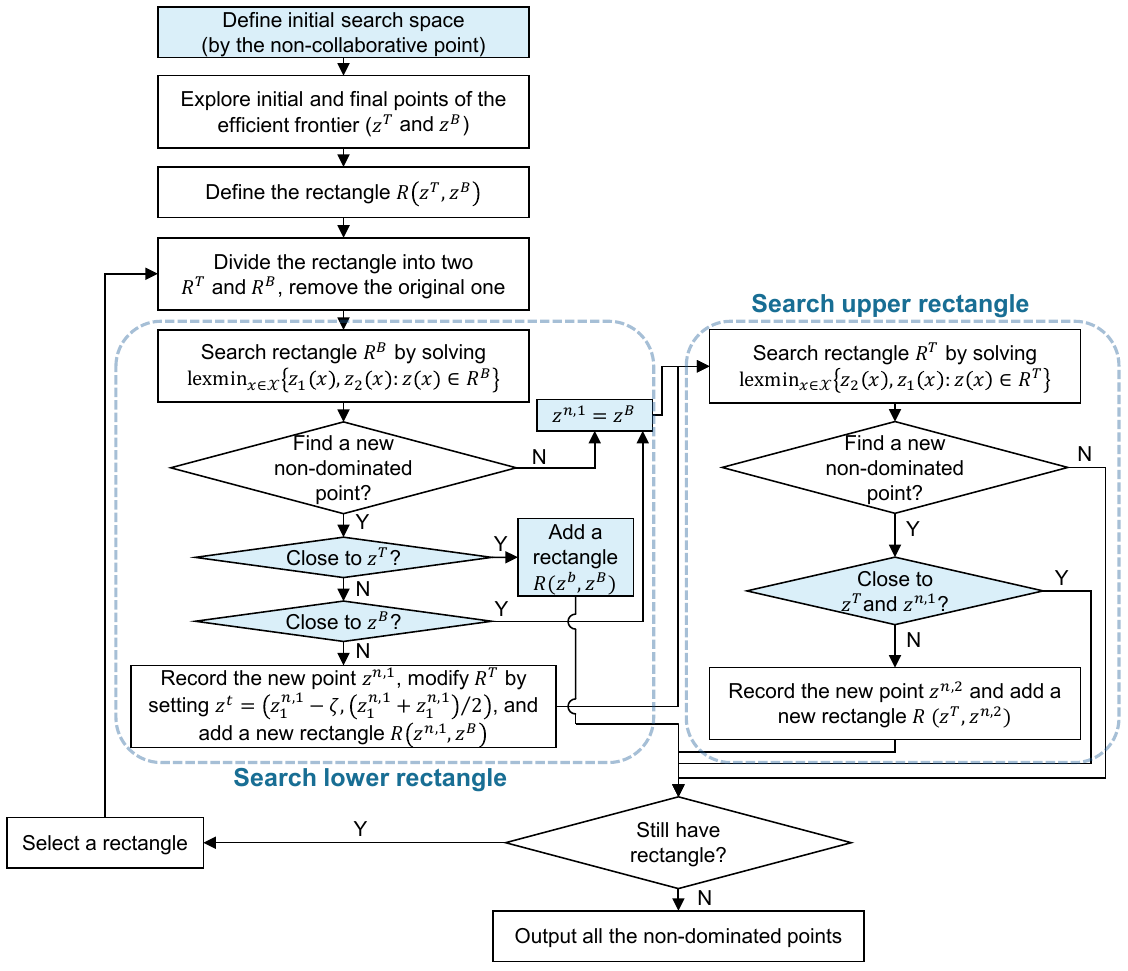}
    \caption{The flowchart of B3M1}
	\label{FIG:flowchart1}
\end{figure*}

The flowchart in Fig. \ref{FIG:flowchart1} illustrates the detailed procedures for B3M1. The steps that have been revised or newly introduced are highlighted in blue, providing a clear distinction between the original and modified processes. As shown, the most significant modification is the addition of judgment steps for evaluating closeness based on strict conditions.

\subsubsection{Balanced bounding box method 2}

\textit{The Balanced Bounding Box Method 2 (B3M2) is based on the relaxed closeness condition defined in Definition 2.2, where solutions satisfying the condition $z^{\rm n}\sim \mathcal{Z}^{\rm e}$ are disregarded. Similar to the original balanced box method, this method requires the definition of rectangles. However, unlike the original method, where the rectangle is defined by $z^T$ and $z^B$ (as shown in Fig. \ref{FIG:1} (a)), the rectangle in this method is defined by the relaxed boundaries $\tilde{z^T} = (z_1^T+\sigma_1, z_2^T-\sigma_2)$ and $\tilde{z^B} = (z_1^B-\sigma_1, z_2^B+\sigma_2)$, as shown in Fig. \ref{FIG:1} (c). The method explores the middle region between $\tilde{z^T}$ and $\tilde{z^B}$, ignoring solutions that satisfy the relaxed closeness condition.}

\begin{figure*} [htbp]
	\centering
		\includegraphics[scale=.6]{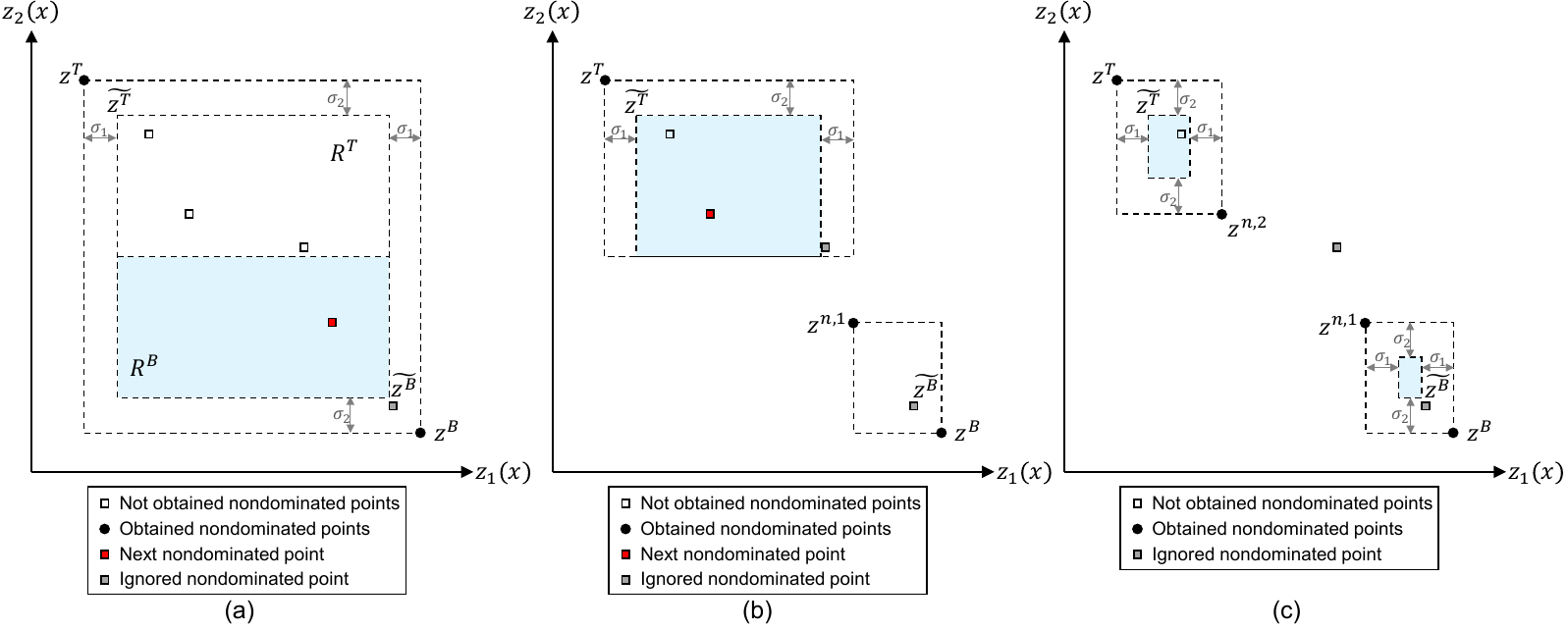}
    \caption{B3M2, $z^{\rm n}\sim \mathcal{Z}^{\rm e}$, (a) investigating $R^B$, (b) investigating $R^T$, (c) new rectangles}
	\label{FIG:M2_1}
\end{figure*}

The procedures of the B3M2, based on the relaxed closeness condition, are depicted in Fig. \ref{FIG:M2_1}.

As depicted in Fig. \ref{FIG:M2_1}, the inner blue rectangle is smaller than the outer transparent rectangle, indicating a considerable reduction of the search space at each step in the B3M compared to the original one. The lower rectangle $\tilde{R^B}$ is defined by the points $\tilde{z^b}$ and $\tilde{z^B}$, where $\tilde{z^b}$ is $(\tilde{z^T_1},(\tilde{z^T_2}+\tilde{z^B_2})/2)$. Upon searching the lower rectangle $\tilde{R^B}$, the red point in Fig. \ref{FIG:M2_1}(a) is obtained, denoted as $z^{\rm n,1}$. Note that if a new non-dominated point $z^{\rm n,1}$ is discovered in the lower rectangle, the upper rectangle $\tilde{R^T}$ should be adjusted in relation to the newly obtained point. The upper rectangle $\tilde{R^T}(\tilde{z^T},\tilde{z^t})$ should be modified by setting $\tilde{z^t}=\left(z^{\rm n,1}_1-\sigma_1,{\rm max}\left(z^{\rm n,1}_2 + \sigma_2, (\tilde{z^T_2}+\tilde{z^B_2})/2\right)\right)$. As illustrated in Fig. \ref{FIG:M2_1}(b), if $z^{\rm n,1}_2 + \sigma_2 \leq (\tilde{z^T_2}+\tilde{z^B_2})/2$, then the lower bound for objective 2 is $(\tilde{z^T_2}+\tilde{z^B_2})/2$. The search area in the upper rectangle is represented by the blue area, with the red point indicating the newly obtained solution. Furthermore, the gray point in the lower rectangle of  Fig. \ref{FIG:M2_1}(b) is excluded due to its proximity to the non-dominated point $z_B$. Similarly, in Fig. \ref{FIG:M2_1}(c), the gray point in the middle is excluded because of its closeness to the non-dominated point $z^{\rm n,1}$. This results in the identification of two new rectangles, shown as the blue areas in Fig. \ref{FIG:M2_1} (c). In this case, the B3M yields two fewer non-dominated points compared to the original method.

The flowchart illustrating the detailed procedures for B3M2 is presented in Fig. \ref{FIG:flowchart}. The purple highlights clearly distinguish the revised or newly introduced steps from the original processes. As shown, the most significant modification involves redefining the rectangles.

\begin{figure*} [htbp]
	\centering
		\includegraphics[scale=0.58]{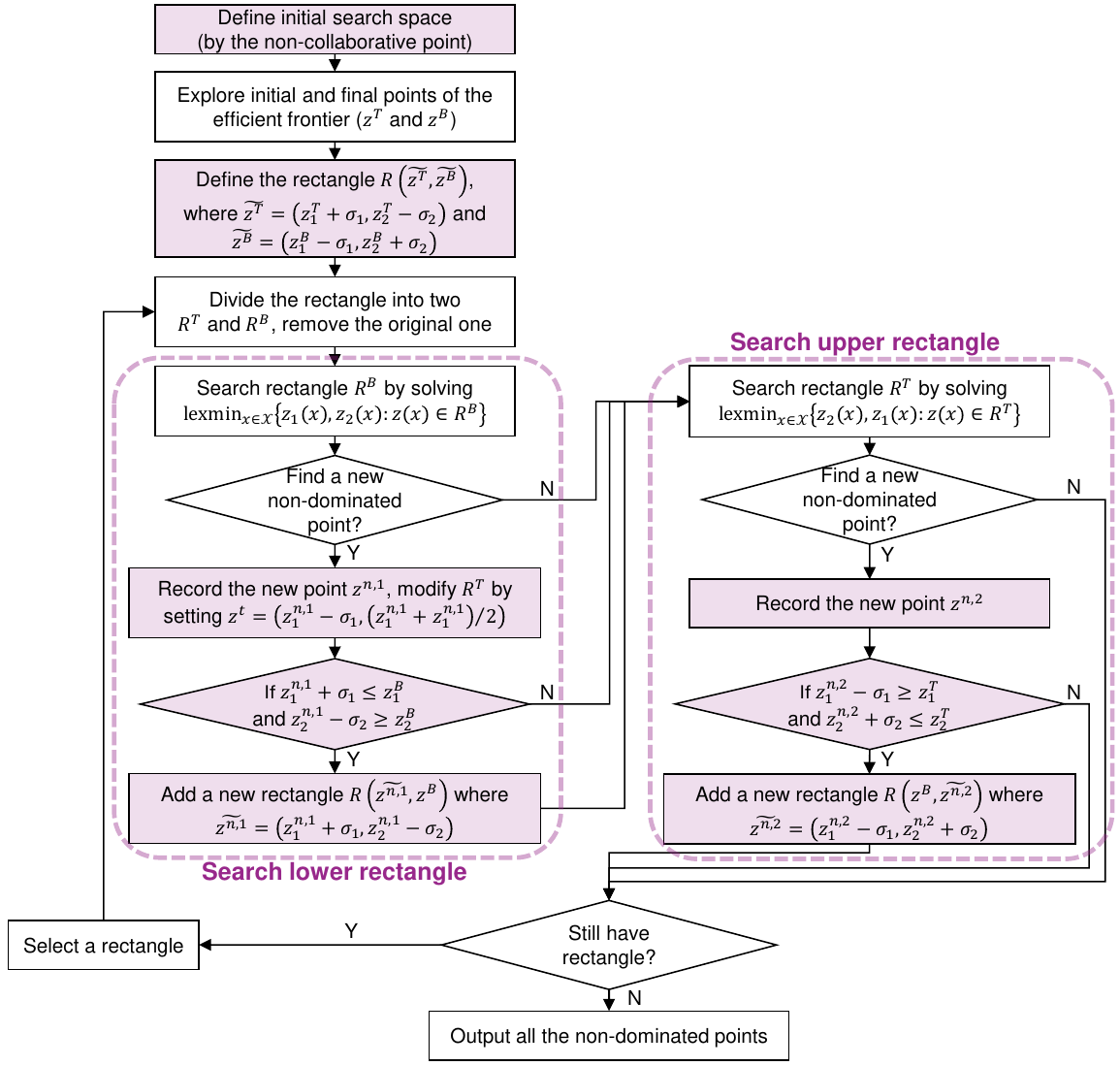}
    \caption{The flowchart of B3M2}
	\label{FIG:flowchart2}
\end{figure*}

\smallskip

\noindent \textbf{\textit{Definition 3: Tolerance range}}

\textit{The term "tolerance range," as represented by the variable $\epsilon$, expressed as a percentage, is that which is required to ensure consistency in the scaling of two objectives. The optimal solutions for these objectives are represented by $z_1^T$ and $z_2^B$, respectively. In consequence, the tolerance margins for objectives 1 and 2 are defined as follows: $\sigma_1=\epsilon z_1^T$ and $\sigma_2=\epsilon z_2^B$. These margins serve to establish the permissible deviation from the optimal solutions, thereby quantifying the degree of acceptance within the specified tolerance range.}

\subsection{Cooperative bargaining}

The solution approaches proposed in Section \ref{Section3.1} can derive an efficient frontier. This section is dedicated to the identification of a final solution from the obtained frontier. Given that the two objectives in the bi-objective optimization are neither fully aligned nor entirely opposed, it is possible for players to collaborate in order to reach a mutually agreeable solution. The process is facilitated by the introduction of cooperative bargaining. We use the Nash bargaining solution (NBS), which provides a unique outcome in such games, to ensure an enforceable agreement based on specific axioms \citep{nash1953}.

Let $\mathscr{x}^*$ represent the final solution selected from the set of all possible solutions $\mathscr{X}$ or the set of reduced solutions $\tilde{\mathscr{X}}$ through cooperative bargaining. The corresponding objective function value, $z^* = z(\mathscr{x}^*)$, is denoted in the objective space $\mathcal{Z}$ or reduced objective space $\tilde{\mathcal{Z}}$.

\begin{figure*} [htbp]
	\centering
		\includegraphics[scale=.65]{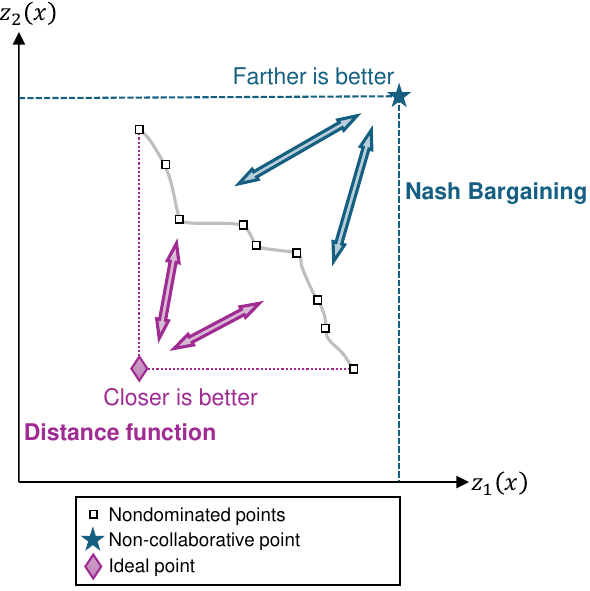}
    \caption{Illustration of cooperative bargaining}
	\label{FIG:bargaining}
\end{figure*}

Before proceeding with a detailed examination of the cooperative bargaining techniques, it is crucial to underscore the fundamental objective of such bargaining: to identify a final solution that represents the greatest distance from the non-collaborative point while simultaneously approaching the ideal point. The aforementioned concepts are illustrated in Fig. \ref{FIG:bargaining} and correspond to the two bargaining methods, namely generalized Nash bargaining and the distance function.

\subsubsection{Generalized Nash bargaining}

The initial step is to compare each solution on the frontier with the non-collaborative solution. This allows us to identify the final agreement point. The optimal final agreement point is the one farthest from the non-collaborative solution.
In scenarios where the two players possess unequal bargaining power, let $\pi$ and $1 - \pi$ denote the relative bargaining strengths of Player 1 and Player 2, respectively. Under these conditions, the Nash bargaining solution is given by the following formulation:

\vspace{-\topsep}
\begin{equation} \label{GNB}
{\rm s}^*= \mathop{{\rm argmax}}\limits_{\left(z_1,z_2\right)\in \mathcal{Z}} (\overline{z}_1 - z_1)^{\pi}(\overline{z}_2 - z_2)^{1-\pi},
\end{equation}

\noindent where $\overline{z}_1$ and $\overline{z}_2$ denote the costs incurred by both companies in the non-cooperative scenario, indicating the maximum cost that they would be willing to accept for such collaboration.
In the symmetric case, the exponents are both equal to 1/2 and sum to 1, thereby ensuring normalization. These exponents represent the bargaining power of each player. As $\pi$ increases, player 1 receives a larger share in equilibrium, and vice versa.

\subsubsection{Distance function}

An alternative method for identifying the final agreement point is to evaluate the solutions on the frontier in comparison to the ideal solution. The optimal final agreement point is thus defined as the solution that is closest to the ideal solution. A distance function is applied to quantify the proximity between a solution on the frontier and the ideal solution, as demonstrated by \citet{zhong2020}, and is defined as follows:

\vspace{-\topsep}
\begin{equation} \label{DisFun}
	D_\alpha \left(z_1,z_2\right)=\left[\left(\frac{z_1-\underline{z}_1}{\overline{z}_1-\underline{z}_1}\right)^{\alpha}+\left(\frac{z_2-\underline{z}_2}{\overline{z}_2-\underline{z}_2}\right)^{\alpha}\right]^{1/\alpha},
\end{equation}

\noindent where $\underline{z}_1$ and $\underline{z}_2$ represent the minimal costs for the players, which are determined by individually minimizing each objective. The parameter $\alpha$ represents the distance parameter. In Eq. \eqref{DisFun}, a min-max normalization is applied to rescale the values.

The $p$-norm (also called $l_p$-norm) of vector $f(z)=(f(z_1),...,f(z_n))$ is:

\vspace{-\topsep}
\begin{equation} \label{DefinePnorm}
    \Vert{f(z)}\Vert_p:=\left(\sum_{i=1}^n\left|f(z_i)\right|^p\right)^{1/p},
\end{equation}

\noindent where the functions $f(z_1)$ and $f(z_2)$ are used to represent the expressions $\frac{z_1-\underline{z}_1}{\overline{z}_1-\underline{z}_1}$ and $\frac{z_2-\underline{z}_2}{\overline{z}_2-\underline{z}_2}$, respectively. The variable $\alpha$ is used to represent the parameter $p$. In practice, $\alpha$ can take any positive real value depending on the desired sensitivity or distance measure.
Given a specific value of $\alpha$, the point that is closest to the ideal point, which represents the Nash bargaining solution, can be determined by the following equation:

\vspace{-\topsep}
\begin{equation} \label{DisFun1}
	{\rm{s}}^*=\mathop{\rm{argmin}}\limits_{z\in \mathcal{Z}} \Vert{f(z)}\Vert_\alpha:=\mathop{\rm{argmin}}\limits_{\left(z_1,z_2\right)\in \mathcal{Z}} D_\alpha \left(z_1,z_2\right).
\end{equation}

For $\alpha=1$, the metric reduces to the taxicab norm (or Manhattan norm), while for $\alpha=2$, it corresponds to the Euclidean norm, representing the straight-line distance from the origin to the point.

To analyze the general properties of the 
$\alpha$-norm, we derive its upper and lower bounds based on the parameter $\alpha$. Specifically, the bounds for $\Vert{f(z)}\Vert_\alpha$ are given by:

\vspace{-\topsep}
\begin{equation} \label{UpperPnorm}
    \left(\sum_{i=1}^n\left|f(z_i)\right|^\alpha\right)^{1/\alpha} \leq \left(\sum_{i=1}^n\left|\mathop{\rm{max}}\limits_i f(z_i)\right|^\alpha\right)^{1/\alpha} = n^{1/\alpha}\mathop{\rm{max}}\limits_i\left|f(z_i)\right|,
\end{equation}

\noindent and

\vspace{-\topsep}
\begin{equation} \label{LowerPnorm}
    \left(\sum_{i=1}^n\left|f(z_i)\right|^\alpha\right)^{1/\alpha} \geq \left(\mathop{\rm{max}}\limits_i \left|f(z_i)\right|^\alpha\right)^{1/\alpha} = \mathop{\rm{max}}\limits_i\left|f(z_i)\right|.
\end{equation}

These bounds establish the range of the $\alpha$-norm, illustrating its confinement between the maximum individual value and a scaled aggregation determined by $\alpha$. According to Eqs. \eqref{UpperPnorm} and \eqref{LowerPnorm}, $\Vert{f(z)}\Vert_\alpha$ limits by:

\vspace{-\topsep}
\begin{equation} \label{ULPnorm}
    \mathop{\rm{max}}\limits_i\left|f(z_i)\right| \leq \Vert{f(z)}\Vert_\alpha \leq n^{1/\alpha}\mathop{\rm{max}}\limits_i\left|f(z_i)\right|.
\end{equation}

As $\alpha$ approaches $\infty$, the $p$-norm approaches the infinity norm or maximum norm. This behavior highlights how the $\infty$-norm emphasizes the largest deviation, making it particularly suitable for robust optimization problems.

\vspace{-\topsep}
\begin{equation} \label{PnormInfinity}
    \Vert{f(z)}\Vert_\infty:=\mathop{\rm{max}}\limits_i\left|f(z_i)\right| = \mathop{\rm{max}}\limits_i \left|\frac{z_i-\underline{z}_i}{\overline{z}_i-\underline{z}_i}\right|,
\end{equation}

\noindent the corresponding Nash bargaining solution is: 

\vspace{-\topsep}
\begin{equation} \label{DisFun2}
	{\rm{s}}^*=\mathop{\rm{argmin}}\limits_{z\in \mathcal{Z}}\Vert{f(z)}\Vert_\infty:=\mathop{\rm{argmin}}\limits_{z_i\in \mathcal{Z}}
 \left(\mathop{\rm{max}}\limits_i \left|\frac{z_i-\underline{z}_i}{\overline{z}_i-\underline{z}_i}\right|\right).
\end{equation}

Based on the properties of the $\infty$-norm, Eq. \eqref{DisFun2} formulates a min-max optimization problem, which is commonly used to handle worst-case scenarios. This approach is commonly used in equilibrium problems and robust optimization. Overall, the $\alpha$-norm framework provides a flexible framework to measure distances, ranging from the Manhattan norm ($\alpha=1$) to the Euclidean norm ($\alpha=2$), and ultimately to the maximum norm ($\alpha\to\infty$), each offering unique properties for different sensitivities in optimization problems.

\section{Case Study} \label{Section4}

In this section, we demonstrate the proposed model and solution methods using case studies based on real-world scenarios. 

\subsection{Scenario setting}

\begin{figure*} [htbp]
	\centering
		\includegraphics[scale=.18]{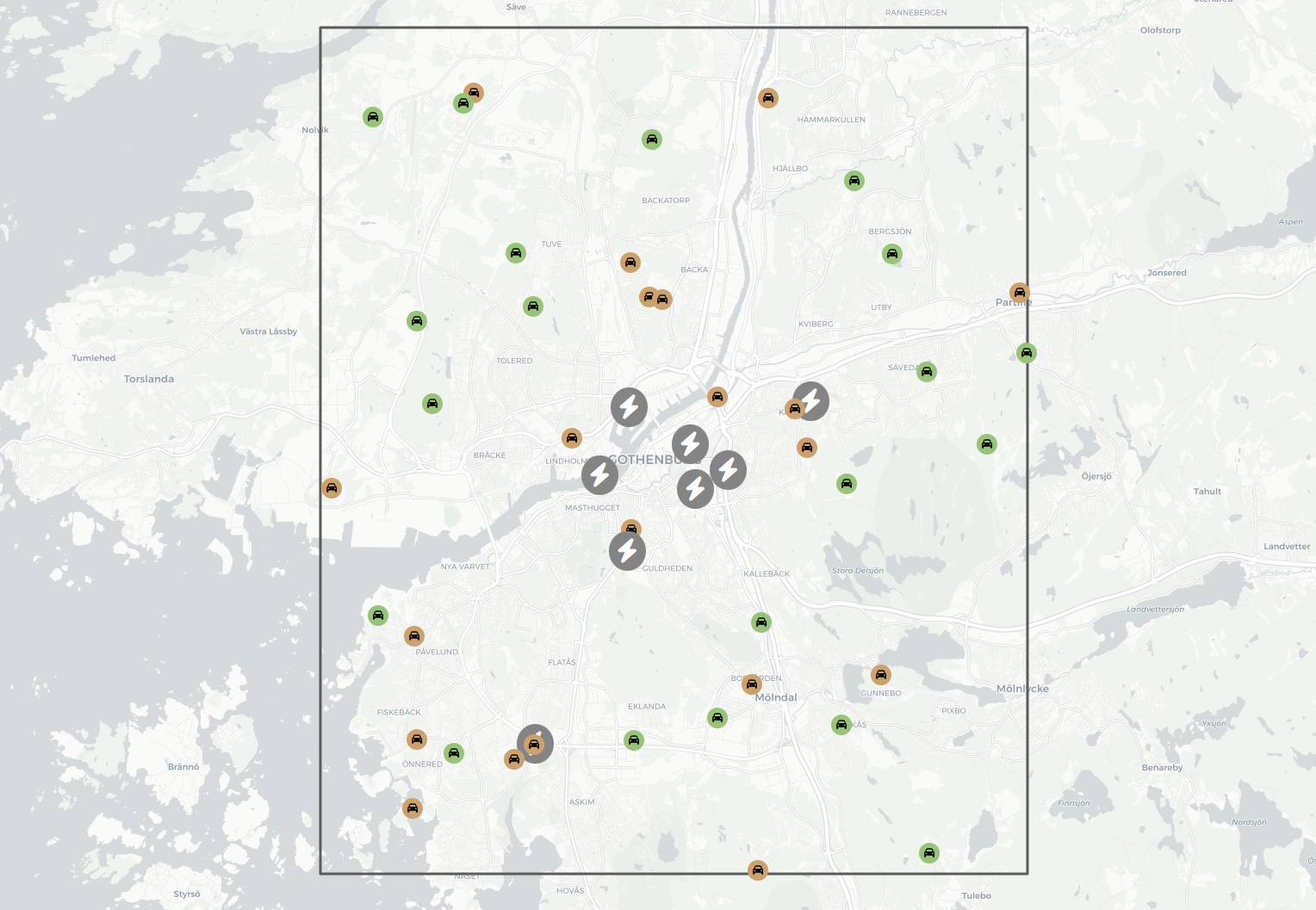}
    \caption{Geographical distribution of charging stations and EVs}
	\label{FIG:map}
\end{figure*}

The case studies utilize real data from charging stations operated by Göteborg Energi, the largest energy provider in Gothenburg, Sweden. In this study, Göteborg Energi acts as a third-party provider, while the collaboration takes place between two fleet operators. The dataset includes precise station locations, along with details on charging power and pricing \citep{GoteEnergi2024}. Fig. \ref{FIG:map} illustrates the charging station locations and the randomly distributed EVs within a 20 km × 20 km service area, where the framed region represents the service boundary. The two fleet operators are represented in green and orange, respectively, whereas the charging stations available for rent are shown in grey.
The distance traveled by vehicles to reach these charging stations is calculated using a real-world map and leveraging the OpenRouteService (ORS) API \citep{ORS2024}.

\begin{figure*} [htbp]
	\centering
		\includegraphics[scale=.4]{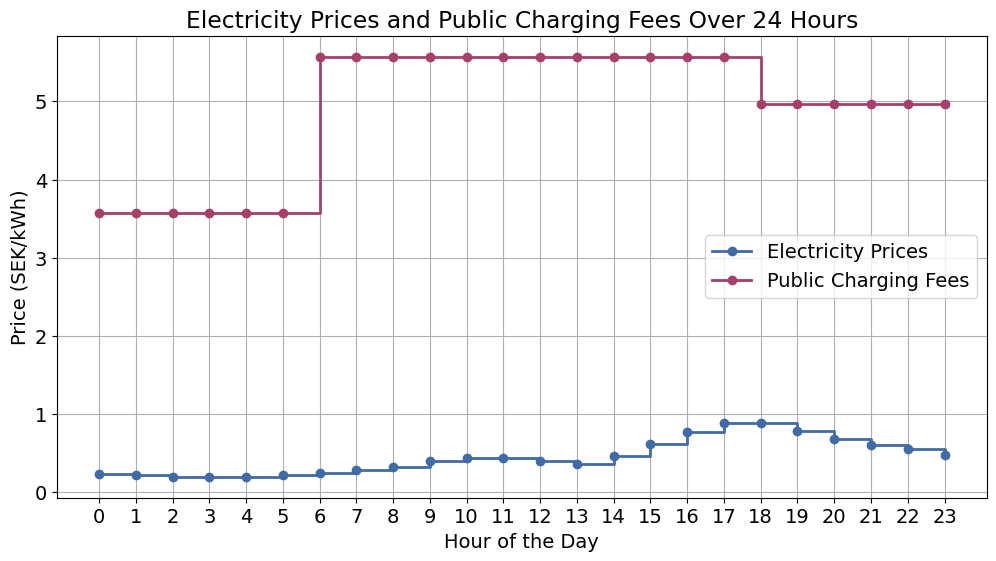}
    \caption{Hourly fluctuations in electricity prices vs. public charging fees}
	\label{FIG:hoursfee}
\end{figure*}

Given the inherent uncertainty in charging demand, we simulate charging requests from EVs distributed across the service area.
The hourly electricity prices are taken from real data provided by Nord Pool \citep{Nordpool2024}, and the public charging fee represents the prices charged by Göteborg Energi's charging stations. Fig. \ref{FIG:hoursfee} illustrates an example of fluctuating electricity prices and public charging fees over a 24-hour period, revealing the temporal price trends at different times of the day.

Two fleet operators collaborate within the system, denoted as $K$. Each fleet operator $k$ manages a group of EVs, $I_k$, while Göteborg Energi provides a set of chargers, $J$. The service period under consideration spans one day, which is divided into 24 one-hour intervals ($T$) to allow for precise time-specific scheduling adjustments. Each charging station incurs a fixed daily rental fee, ${\rm f}_j^k$, which is calculated on a per-day basis. This reflects a flexible rental strategy that is ideally suited to accommodate fluctuating demand. The fixed daily rental fee is set at 1500 SEK, based on the public charging fee and roughly equivalent to the cost of using a public charging station for 7 hours. To maintain general applicability, this study considers only charging stations from Göteborg Energi with a maximum charging rate of 50 kW \citep{GoteEnergi2024}, which corresponds to the previously defined charging rate $\varepsilon_{ij}$. In scenarios with more stable demand, alternative business models (such as subscription-based models, long-term leasing agreements, or even the construction of new facilities) could be considered. When a company rents a charging station, it is required to pay the electricity fee, denoted as $c_{\rm e,o}^{j,\tau}$. For collaborative EVs, a discount factor $\beta$ of 0.5 is applied, meaning that a collaborating company pays half of the public charging fee when utilizing charging stations rented by others. This corresponds to the previously defined unit energy fee for collaborative EVs, $c_{\rm e,c}^{j,\tau}$. Additionally, the cost of waiting for charging, represented by ${\lambda}_i$, is estimated at 300 SEK per hour. The energy cost for EV $i$ traveling to charger $j$, denoted by ${w}_{ij}$, is calculated based on the unit energy consumption cost of 6 SEK per kilometer.

In order to assess the performance at different scales, three cases with varying numbers of EVs were examined. In Case \#1, there are 20 EVs, in Case \#2, there are 30 EVs, and in Case \#3, there are 40 EVs.

\subsection{Computational performance}

\begin{figure*} [htbp]
	\centering
            \subfigure[]{
		  \includegraphics[width=0.14\textwidth]{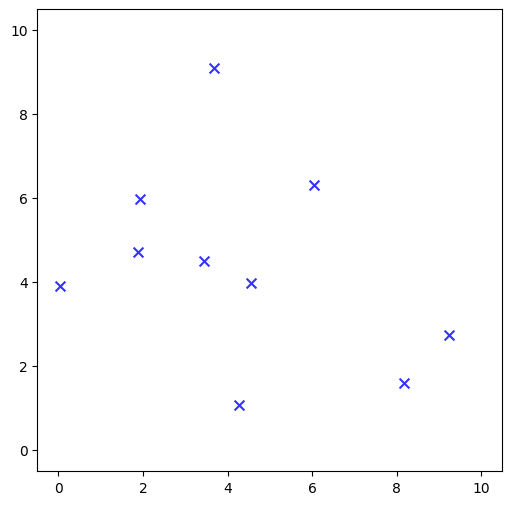}}
            \hspace{0.02in}
            \subfigure[]{
            \includegraphics[width=0.14\textwidth]{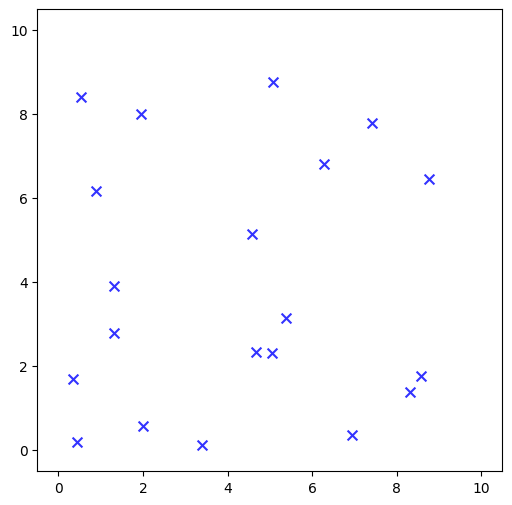}}
            \hspace{0.02in}
            \subfigure[]{
            \includegraphics[width=0.14\textwidth]{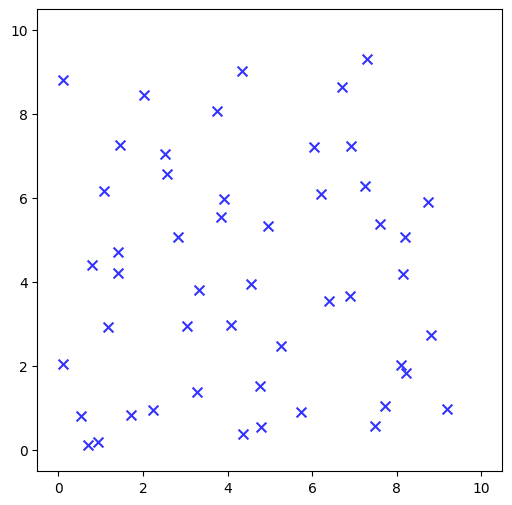}}
            \hspace{0.02in}
            \subfigure[]{
            \includegraphics[width=0.14\textwidth]{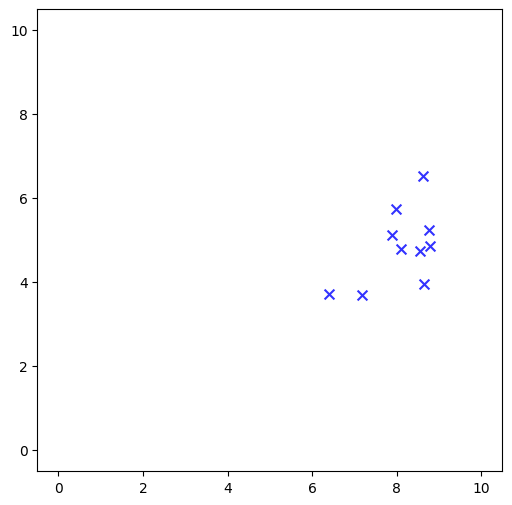}}
            \hspace{0.02in}
            \subfigure[]{
            \includegraphics[width=0.14\textwidth]{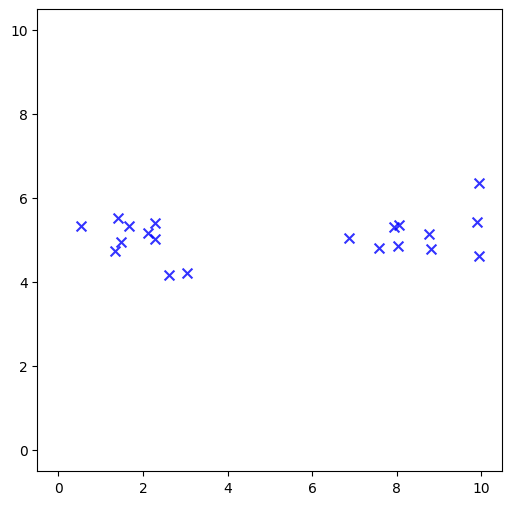}}
            \hspace{0.02in}
            \subfigure[]{
            \includegraphics[width=0.14\textwidth]{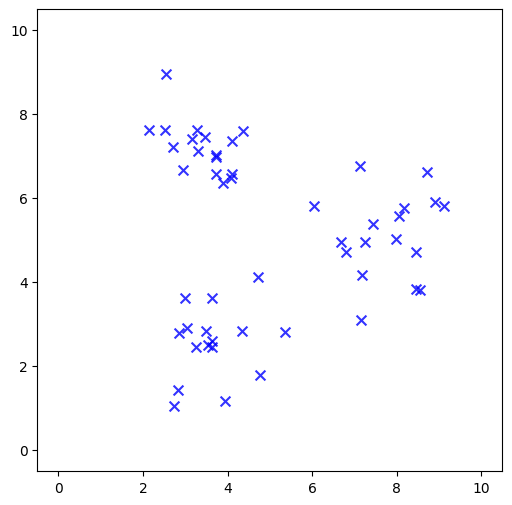}}
            \hspace{0.02in}
            \subfigure[]{
		  \includegraphics[width=0.15\textwidth]{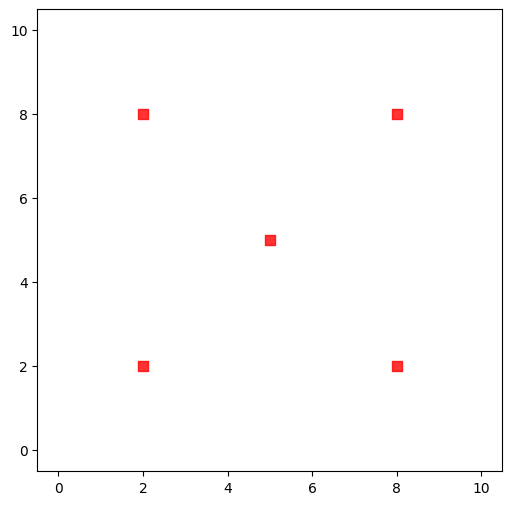}}
            \hspace{0.02in}
            \subfigure[]{
            \includegraphics[width=0.15\textwidth]{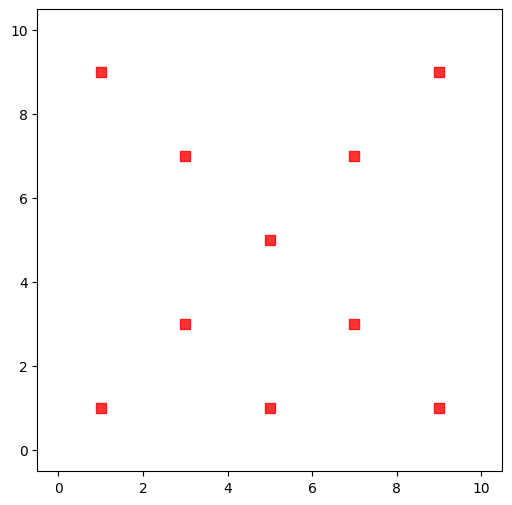}}
            \hspace{0.02in}
            \subfigure[]{
            \includegraphics[width=0.15\textwidth]{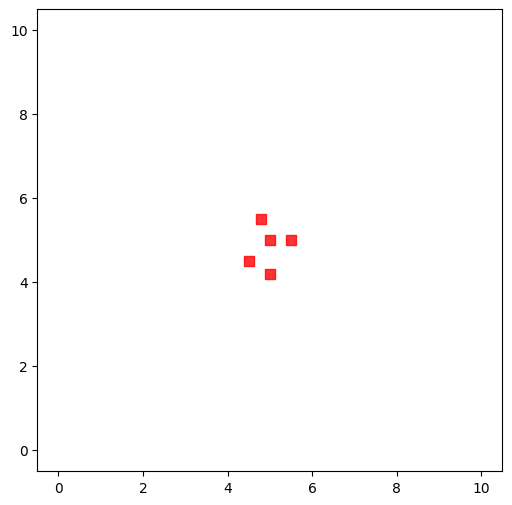}}
            \hspace{0.02in}
            \subfigure[]{
            \includegraphics[width=0.15\textwidth]{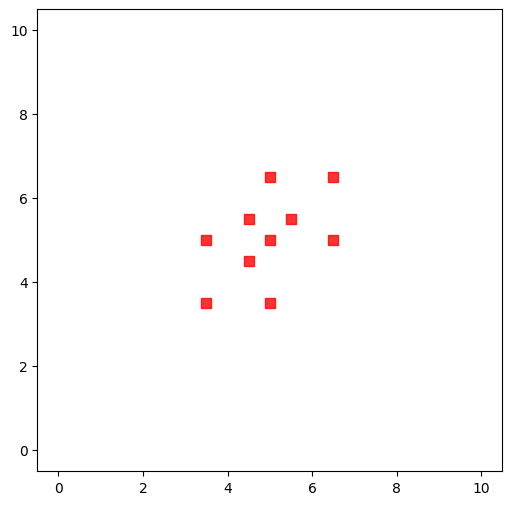}}
        \caption{
    Distribution of EVs and chargers in test cases. 
    Subfigures (a)–(f) show electric vehicle (EV) spatial distributions with different fleet sizes (10, 20, 50) under uniform and clustered configurations. 
    Subfigures (g)–(j) show charger layouts with two configurations (uniform and centralized) and varying charger counts (5 or 10).}
	\label{FIG:DisChar}
\end{figure*}

This section presents the computational performance of the developed B3Ms on a set of test cases. Various distributions of EVs and chargers are considered. The EV locations follow two spatial distributions: uniform and clustered, while chargers are arranged using two distinct layouts: uniform and centralized (see Fig. \ref{FIG:DisChar}). Each EV is assigned a time window within a 24-hour period, with the earliest start time randomly generated between hours 1 and 20. The time window length is fixed at 4 hours, and each hour represents one time interval.

To distinguish between different test cases, we adopt the naming convention DistEV–DistChar–NumEV–NumCharger, where DistEV refers to the EV distribution (uniform or clustered), DistChar refers to the charger layout (uniform or centralized), and NumEV and NumCharger indicate the number of EVs and chargers, respectively. For example, UniEV–UniChar–50–10 represents a case where EVs are uniformly distributed, chargers follow a uniform layout, and there are 50 EVs and 10 chargers.

We evaluate the performance of B3M1 and B3M2 compared to the original balanced box method under two service area sizes: 10 km × 10 km and 100 km × 100 km. The tolerance range $\epsilon$ is set to 3\%. Key performance metrics include the number of non-dominated points (NDP), computational time (CPU), solution gap percentage (Gap), and computational time savings (CTS). The Gap metric quantifies the quality loss of B3M-generated solutions compared to the original method. For each ignored solution (i.e., one not retained by B3M), we compute the Euclidean distance in the objective space to all retained solutions, based on their objective values. The smallest of these distances is recorded as the individual gap. The overall gap is then defined as the average of these individual gaps across all ignored solutions, providing a quantitative measure of the quality loss introduced by filtering out nearby solutions. CTS is calculated as the percentage reduction in CPU time compared to the original method.

As shown in Tables \ref{tbcompresults} and \ref{tbcompresults2}, both B3M1 and B3M2 significantly reduce computational time while maintaining solution quality. B3M2 consistently achieves the highest CTS, especially in larger instances, frequently exceeding 60\% and reaching up to 93.7\%. B3M1 also yields substantial reductions. Both methods maintain comparable or lower NDP than the original method, indicating their ability to eliminate closely positioned solutions without losing diversity. As shown in Table \ref{tbcompresults}, gap values remain within the 3\% tolerance in most 10 km × 10 km cases, with only one exception. However, under the 100 km × 100 km setting (Table \ref{tbcompresults2}), B3M2 exhibits gap values exceeding the threshold in 8 out of 24 test cases, particularly those with clustered EV distributions and centralized charger layouts. This indicates that spatially dispersed scenarios introduce greater solution variability. Nonetheless, the gap values remain reasonably close to 3\%, suggesting acceptable accuracy for large-scale applications. Overall, B3M2 offers the best balance between efficiency and accuracy, demonstrating strong scalability across service areas.

Comparing different layout and distribution settings further reveals that charger layout has a more decisive impact on performance. Centralized layouts concentrate chargers in compact regions, reducing solution space and CPU time, but limiting assignment flexibility. In contrast, uniform layouts increase spatial diversity, leading to more non-dominated solutions and a higher computational burden. EV distribution has less consistent influence across cases. These results underscore the importance of charger placement in scheduling outcomes and validate the robustness of B3M2 under varied conditions. The next section further illustrates the solution characteristics through representative case studies.


\begin{table}[width=.9\linewidth,cols=11,pos=htbp]
\caption{Computational results for the test cases (Service area 10 km × 10 km)}\label{tbcompresults}
\begin{tabular*}{\tblwidth}{@{} L|LL|LLLL|LLLL|@{} }
\toprule
& \multicolumn{2}{c|}{Balanced box} & \multicolumn{4}{c|}{B3M1} & \multicolumn{4}{c|}{B3M2}\\
Test cases & NDP & CPU (s) & NDP & CPU (s) & Gap & CTS & NDP & CPU (s) & Gap & CTS\\
\midrule
UniEV-UniChar-10-5 & 7 & 1314.6 & 5 & 1327.6 & 2.25\% & -1.0\% & 3 & 567.6 & 3.18\% & 56.9\% \\
UniEV-UniChar-10-10 & 12 & 4492.7 & 5 & 2433.5 & 1.30\% & 45.8\% & 3 & 1083.4 & 2.23\% & 75.8\% \\
UniEV-UniChar-20-5 & 4 & 1513.7 & 3 & 1512.4 & 0.84\% & 0.1\% & 2 & 756.4 & 2.80\% & 50.0\% \\
UniEV-UniChar-20-10 &  7 & 5515.9 & 2 & 2037.0 & 2.03\% & 63.1\% & 2 & 1496.1 & 2.03\% & 72.9\%\\
UniEV-UniChar-50-5 & 30 & 31510.8 & 3 & 6996.4 & 1.32\% & 77.8\% & 2 & 1996.6 & 1.97\% & 93.7\% \\
UniEV-UniChar-50-10 & 3 & 10959.7 & 2 & 7636.6 & 0.63\% & 30.3\% & 2 & 4252.3 & 0.63\% & 61.2\% \\
Average & 10.5 & 9217.9 & 3.3 & 3657.3 & 1.40\% & 36.0\% & 2.3 & 1692.6 & 2.14\% & 68.4\% \\ 
\midrule
UniEV-CenChar-10-5 & 7 & 1306.4 & 3 & 623.4 & 0.55\% & 52.3\% & 3 & 571.3 & 0.55\% & 56.6\% \\
UniEV-CenChar-10-10 & 10 & 3679.4 & 3 & 1287.4 & 1.52\% & 65.0\% & 3 & 1152.2 & 1.52\% & 68.7\% \\
UniEV-CenChar-20-5 & 5 & 1905.9 & 2 & 980.6 & 0.65\% & 48.6\% & 2 & 737.4 & 0.65\% & 61.3\% \\
UniEV-CenChar-20-10 & 12 & 94558.3 & 2 & 2161.8 & 0.95\% & 77.1\% & 2 & 1572.4 & 0.95\% & 83.4\%\\
UniEV-CenChar-50-5 & 13 & 13813.7 &  2 & 2373.7 & 0.37\% & 82.8\% & 2 & 1914.7 & 0.37\% & 86.1\% \\
UniEV-CenChar-50-10 & 5 & 20122.8 & 2 & 6155.0 & 0.32\% & 69.4\% & 2 & 4413.8 & 0.32\% & 78.7\% \\
Average & 8.7 & 8381.1 & 2.3 & 2263.6 & 0.73\% & 65.9\% & 2.3 & 1726.4 & 0.73\% & 72.4\% \\ 
\midrule
CluEV-UniChar-10-5 & 4 & 844.5 & 3 & 639.7 & 0.36\% & 24.2\% & 3 & 571.5 & 0.36\% & 32.3\% \\
CluEV-UniChar-10-10 & 4 & 1598.8 & 3 & 1285.9 & 0.46\% & 19.6\% & 3 & 1163.9 & 0.46\% & 27.2\% \\
CluEV-UniChar-20-5 & 11 & 4073.0 & 3 & 1542.5 & 1.40\% & 62.1\% & 3 & 1148.8 & 1.50\% & 71.8\% \\
CluEV-UniChar-20-10 & 10 & 8005.9 & 2 & 2129.9 & 1.17\% & 73.4\% & 2 & 1585.6 & 1.17\% & 80.2\% \\
CluEV-UniChar-50-5 & 8 & 8217.6 & 3 & 5031.5 & 0.86\% & 38.8\% & 2 & 2021.9 & 2.21\% & 75.4\% \\
CluEV-UniChar-50-10 & 7 & 21722.6 & 2 & 9648.3 & 0.96\% & 55.6\% & 2 & 4368.7 & 0.96\% & 79.9\% \\
Average & 7.3 & 7410.4 & 2.7 & 3379.7 & 0.87\% & 45.6\% & 2.5 & 1810.1 & 1.11\% & 61.1\% \\ 
\midrule
CluEV-CenChar-10-5 & 6 & 1136.2 & 3 & 879.9 & 0.37\% & 22.6\% & 3 & 580.2 & 0.37\% & 48.9\%\\
CluEV-CenChar-10-10 & 4 & 1468.1 & 3 & 1277.0 & 0.47\% & 13.0\% & 3 & 1139.8 & 0.47\% & 22.4\% \\
CluEV-CenChar-20-5 & 4 & 1476.6 & 2 & 896.4 & 0.19\% & 39.3\% & 2 & 766.6 & 0.19\% & 48.1\%\\
CluEV-CenChar-20-10 & 3 & 2367.9 & 2 & 1924.9 & 0.48\% & 18.7\% & 2 & 1569.7 & 0.48\% & 33.7\%\\
CluEV-CenChar-50-5 & 5 & 5422.7 & 2 & 2845.0 & 0.26\% & 47.5\% & 2 & 2008.6 & 0.26\% & 63.0\%\\
CluEV-CenChar-50-10 & 9 & 30620.0 & 2 & 6116.7 & 0.63\% & 80.0\% & 2 & 4470.0 & 0.63\% & 85.4\%\\
Average & 5.2 & 7081.9 & 2.3 & 2323.3 & 0.40\% & 36.9\% & 2.3 & 1755.8 & 0.40\% & 50.2\% \\ 
\midrule
Total average & 7.9 & 8022.8 & 2.7 & 2905.9 & 0.85\% & 46.1\% & 2.4 & 1746.2 & 1.09\% & 63.1\% \\ 
\bottomrule
\end{tabular*}
\end{table}


\begin{table}[width=.9\linewidth,cols=11,pos=htbp]
\caption{Computational results for the test cases (Service area 100 km × 100 km)}\label{tbcompresults2}
\begin{tabular*}{\tblwidth}{@{} L|LL|LLLL|LLLL|@{} }
\toprule
& \multicolumn{2}{c|}{Balanced box} & \multicolumn{4}{c|}{B3M1} & \multicolumn{4}{c|}{B3M2}\\
Test cases & NDP & CPU (s) & NDP & CPU (s) & Gap & CTS & NDP & CPU (s) & Gap & CTS\\
\midrule
UniEV-UniChar-10-5 & 2 & 379.2 & 2 & 361.6 & 0\% & 4.6\% & 2 & 361.2 & 0\% & 4.7\% \\
UniEV-UniChar-10-10 & 3 & 1122.1 & 2 & 845.3 & 0.59\% & 24.7\% & 1 & 480.1 & 4.68\% & 57.2\% \\
UniEV-UniChar-20-5 & 6 & 2200.2 & 5 & 1959.2 & 1.01\% & 11.0\% & 3 & 1090.6 & 3.84\% & 50.4\% \\
UniEV-UniChar-20-10 & 7 & 5500.7 & 6 & 4772.5 & 1.64\% & 13.2\% & 4 & 2984.9 & 3.47\% & 45.7\% \\
UniEV-UniChar-50-5 & 34 & 34148.2 & 5 & 13123.4 & 1.69\% & 61.6\% & 4 & 3825.9 & 2.52\% & 88.8\% \\
UniEV-UniChar-50-10 & 15 & 31090.3 & 3 & 10995.1 & 1.92\% & 64.6\% & 3 & 6101.0 & 1.92\% & 80.4\% \\
Average & 11.2 & 12406.8 & 3.8 & 5342.9 & 1.14\% & 30.0\% & 2.8 & 2473.9 & 2.74\% & 54.5\% \\ 
\midrule
UniEV-CenChar-10-5 & 6 & 1092.2 & 2 & 655.2 & 2.35\% & 40.0\% & 2 & 360.4 &2.35\% & 67.0\% \\
UniEV-CenChar-10-10 & 4 & 1472.4 & 2 & 1191.5 & 2.69\% & 19.1\% & 2 & 709.6 & 2.69\% & 51.8\%\\
UniEV-CenChar-20-5 & 5 & 1833.3 & 3 & 1586.9 & 2.78\% & 13.4\% & 3 & 1100.3 & 2.78\% & 40.0\%\\
UniEV-CenChar-20-10 & 11 & 8386.9 & 8 & 7630.2 & 1.20\% & 9.0\% & 2 & 1937.6 & 2.51\% & 55.9\%\\
UniEV-CenChar-50-5 & 17 & 17386.2 & 3 & 7035.2 & 1.30\% & 59.5\% & 2 & 1937.6 & 1.57\% & 88.9\%\\
UniEV-CenChar-50-10 & 14 & 30377.8 & 4 & 18781.6 & 1.05\% & 38.2\% & 2 & 4294.1 & 2.94\% & 85.9\%\\
Average & 9.5 & 10091.5 & 3.7 & 6146.8 & 1.90\% & 29.9\% & 2.7 & 2017.1 & 2.47\% & 64.9\% \\ 
\midrule
CluEV-UniChar-10-5 & 3 & 565.4 & 2 & 449.7 & 1.57\% & 20.5\% & 1 & 252.6 & 7.60\% & 55.3\% \\
CluEV-UniChar-10-10 & 3 & 1137.8 & 3 & 1169.6 & 0\% & -2.8\% & 2 & 774.9 & 2.38\% & 31.9\% \\
CluEV-UniChar-20-5 & 5 & 1900.5 & 4 & 1660.3 & 0.39\% & 12.6\% & 3 & 1136.1 & 3.60\% & 40.2\% \\
CluEV-UniChar-20-10 & 6 & 4700.7 & 4 & 4441.7 & 2.28\% & 5.5\% & 3 & 2331.0 & 3.04\% & 50.4\% \\
CluEV-UniChar-50-5 & 21 & 21344.9 & 4 & 9060.6 & 1.77\% & 57.6\% & 3 & 2982.9 & 2.47\% & 86.0\% \\
CluEV-UniChar-50-10 & 33 & 71007.5 & 5 & 29191.0 & 1.48\% & 58.9\% & 4 & 8438.8 & 2.16\% & 88.1\% \\
Average & 11.8 & 16776.1 & 3.7 & 7648.7 & 1.24\% & 26.6\% & 2.7 & 2652.7 & 3.54\% & 58.7\% \\ 
\midrule
CluEV-CenChar-10-5 & 6 & 1186.9 & 4 & 820.0 & 0.83\% & 30.9\% & 3 & 565.2 & 3.24\% & 52.4\% \\
CluEV-CenChar-10-10 & 4 & 1519.6 & 4 & 1451.5 & 0\% & 4.5\% & 3 & 1137.2 & 3.66\% & 25.2\% \\
CluEV-CenChar-20-5 & 6 & 2325.0 & 3 & 1416.1 & 1.30\% & 39.1\% & 3 & 1151.6 & 1.30\% & 50.5\% \\
CluEV-CenChar-20-10 & 7 & 5541.1 & 6 & 5543.7 & 1.12\% & -0.05\% & 3 & 2355.9 & 7.24\% & 57.5\% \\
CluEV-CenChar-50-5 & 6 & 6115.1 & 2 & 2376.6 & 1.31\% & 61.1\% & 2 & 2006.5 & 1.31\% & 67.2\% \\
CluEV-CenChar-50-10 & 21 & 49225.9 & 4 & 25697.2 & 1.81\% & 47.8\% & 3 & 6572.0 & 2.64\% & 86.6\% \\
Average & 8.3 & 10985.6 & 3.8 & 6217.5 & 1.06\% & 30.6\% & 2.8 & 2298.1 & 3.23\% & 56.6\% \\ 
\midrule
Total average & 10.2 & 12565.0 & 3.8 & 6339.0 & 1.34\% & 29.2\% & 2.8 & 2360.4 & 2.99\% & 58.7\% \\ 
\bottomrule
\end{tabular*}
\end{table}

\subsection{Different tolerance range}

\begin{figure*} [htbp]
	\centering
            \subfigure[Case \#1]{
		  \includegraphics[width=0.7\textwidth]{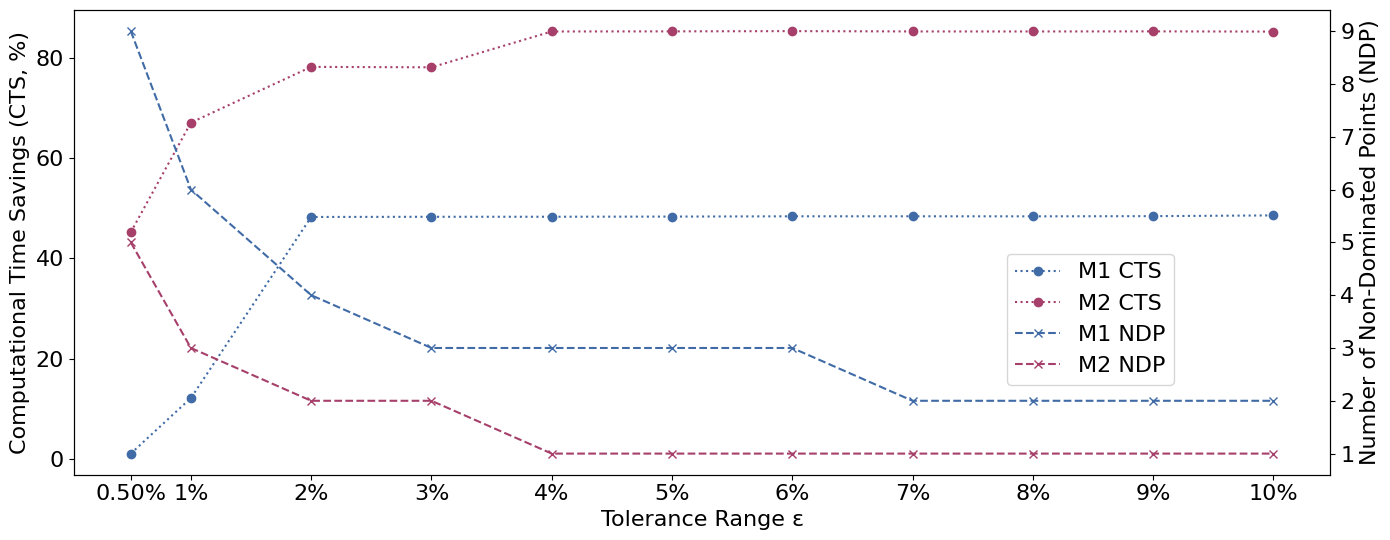}}
            \hspace{0.02in}
            \subfigure[Case \#2]{
		  \includegraphics[width=0.7\textwidth]{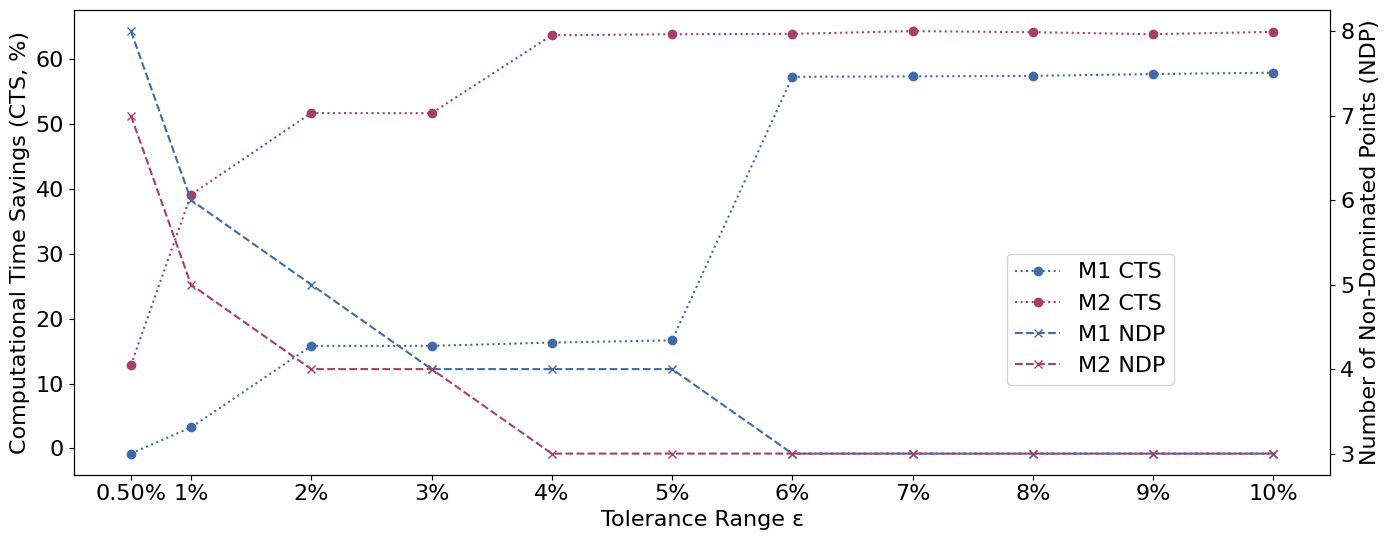}}
            \hspace{0.02in}
            \subfigure[Case \#3]{
            \includegraphics[width=0.7\textwidth]{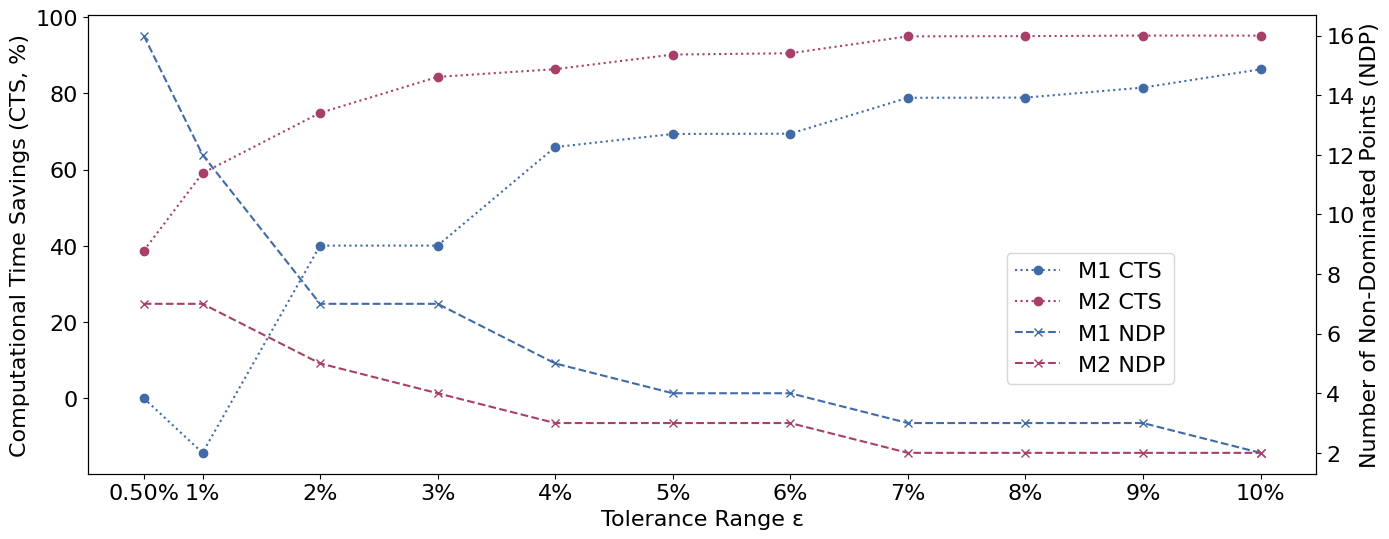}}
        \caption{Comparative analysis of results across various tolerance ranges $\epsilon$}
	\label{FIG:7}
\end{figure*}

This section presents an analysis of the impact of the tolerance range $\epsilon$ on the computational performance and solution sets. Specifically, we examine the three cases defined earlier (Case \#1: 20 EVs, Case \#2: 30 EVs, Case \#3: 40 EVs). By varying the tolerance range $\epsilon$, the relationship between CTS and NDP of three cases is illustrated in Fig. \ref{FIG:7}.

As illustrated in Fig. \ref{FIG:7}, an increase in the tolerance range is generally accompanied by a reduction in the number of non-dominated points and computational time. In all cases, the implementation of B3M2 in comparison to B3M1 has the potential to reduce computational time, although this may be at the cost of a reduction in the number of non-dominated points. The B3M2 has the potential to reduce computational time by over 50\% compared to the original balanced box method when the tolerance range $\epsilon$ is set to at least $2\%$, and up to 95\% when the range is increased further. Notably, when the tolerance range exceeds 6\%, the changes in the solutions are negligible. From an alternative standpoint, the tolerance range represents the level of profit reduction that companies can accept, with reductions exceeding 5\% typically considered significant. Accordingly, the selected tolerance ranges were 1\%, 2\%, and 5\% to facilitate a more comprehensive analysis, as these values provide meaningful variations in results.

Table \ref{tbresults} summarizes the outcomes of various methods, including the balanced box method and its two enhancements: B3M1 and B3M2. The summary provides the number of non-dominated points (NDP), the computational time (CPU), the predetermined tolerance range ($\epsilon$), the gap between the B3M and the original balanced box method (Gap), and the computational time savings (CTS).

\begin{table}[width=.9\linewidth,cols=12,pos=h]
\caption{Results of different methods}\label{tbresults}
\begin{tabular*}{\tblwidth}{@{} L|LLL|LLLL|LLLL|@{} }
\toprule
& \multicolumn{2}{c}{Balanced box method} & \multirow{2}{*}{$\epsilon$} & \multicolumn{4}{c|}{B3M1} & \multicolumn{4}{c|}{B3M2}\\
& NDP & CPU (s) & & NDP & CPU (s) & Gap & CTS & NDP & CPU (s) & Gap & CTS\\
\midrule
\multirow{3}{*}{Case \#1} & \multirow{3}{*}{9} & \multirow{3}{*}{5366.6} & 1\% & 6 & 4778.2 & 0.9\% & 11.0\% & 3 & 1769.2 & 2.0\% & 67.0\% \\
&&& 2\% & 4 & 2778.0 & 1.2\% & 48.2\% & 2 & 1173.5 & 3.8\% & 78.1\%\\
&&& 5\% & 3 & 2774.0 & 2.7\% & 48.3\% & 1 & 793.9 & 7.0\% & 85.2\%\\
\multirow{3}{*}{Case \#2} & \multirow{3}{*}{8} & \multirow{3}{*}{7504.9} & 1\% & 6 & 7261.7 & 1.0\% & 3.2\% &  5 & 4569.2 & 1.4\% & 39.1\% \\
&&& 2\% & 5 & 6319.7 & 1.4\% & 15.8\% & 4 & 3627.1 & 3.4\% & 51.7\%\\
&&& 5\% & 4 & 6255.3 & 3.4\% & 16.7\% & 3 & 2715.9 & 5.2\% & 63.8\%\\
\multirow{3}{*}{Case \#3} & \multirow{3}{*}{17} & \multirow{3}{*}{50644.0} & 1\% & 12 & 57878.2 & 1.1\% & -14.3\% &  7 & 20730.3 & 2.0\% & 59.1\%\\
&&& 2\% & 7 & 30347.3 & 1.7\% & 40.1\% & 4 & 12714.6 & 2.5\% & 74.9\% \\
&&& 5\% & 4 & 15515.6 & 3.0\% & 69.4\% & 3 & 4964.2 & 4.3\% & 90.2\%\\
\bottomrule
\end{tabular*}
\end{table}

As can be seen in Table \ref{tbresults}, the results align with the expectations regarding the design of the B3Ms. As defined in Section \ref{section3.1.2}, the two B3Ms use different criteria for ignoring solutions based on the closeness measures $z^{\rm n} \approx \mathcal{Z}^{\rm e}$ and $z^{\rm n} \sim \mathcal{Z}^{\rm e}$. For B3M1, ignored solutions must have both objectives within the tolerance range $\epsilon$. In contrast, for B3M2, only one objective has to be within the tolerance range $\epsilon$. Several insights can be provided from the results. Firstly, it shows that the number of non-dominated points obtained from B3M1 is greater than that obtained from B3M2.
Secondly, the gaps for B3M1 are generally smaller than the predetermined tolerance range $\epsilon$, with exceptions observed in Case \#3 (1.1\% when $\epsilon=1\%$). In contrast, the gaps for B3M2 frequently exceed the predetermined tolerance range, reflecting greater variability and reduced adherence to $\epsilon$. For example, when $\epsilon=5\%$, the gaps for B3M2 range from 4.3\% (Case \#3) to 7.0\% (Case \#1).
Thirdly, the computational time for both B3Ms is significantly reduced with a larger tolerance range $\epsilon$, with B3M2 demonstrating even more notable savings compared to B3M1. 
However, there is an exception in case \#3 with $\epsilon=1\%$: the computational time of B3M1 is longer than that of the balanced box method. This discrepancy arises because B3M1 disregards the solution by comparing it with the upper-left or lower-right points, which should be recorded as solutions. As a result, the search area is expanded in this particular instance.

\subsection{Distribution of non-dominated points}

\begin{figure*} [htbp]
	\centering
            \subfigure[Complete solution set]{
		  \includegraphics[width=0.23\textwidth]{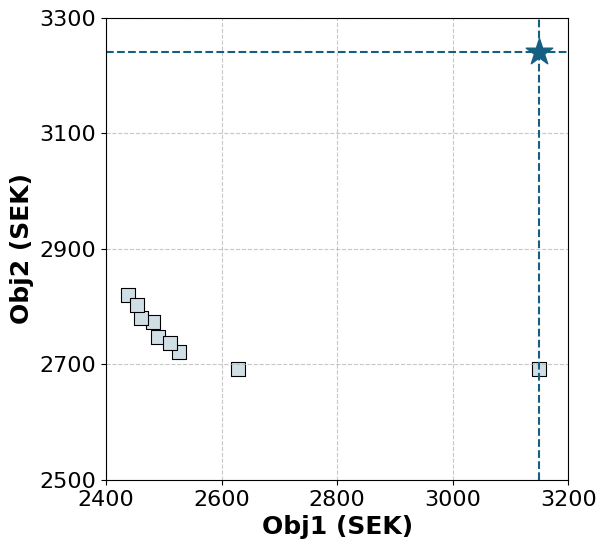}}
            \hspace{0.02in}
            \subfigure[B3M1 with $\epsilon=1\%$]{
            \includegraphics[width=0.23\textwidth]{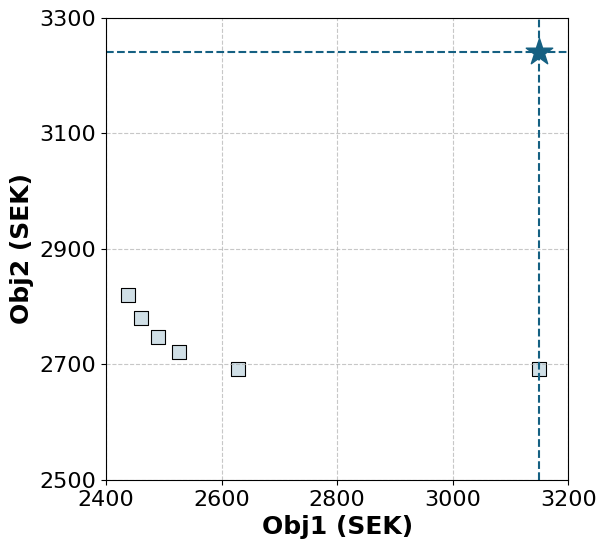}}
            \hspace{0.02in}
            \subfigure[B3M1 with $\epsilon=2\%$]{
            \includegraphics[width=0.23\textwidth]{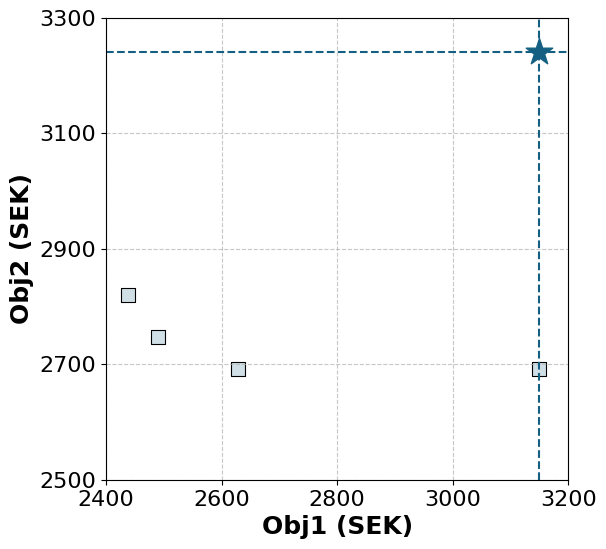}}
            \hspace{0.02in}
            \subfigure[B3M1 with $\epsilon=5\%$]{
            \includegraphics[width=0.23\textwidth]{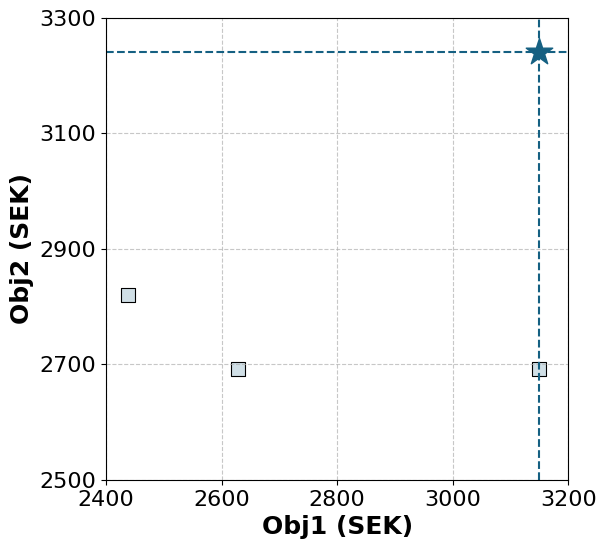}}
            \hspace{0.02in}
            \subfigure[B3M2 with $\epsilon=1\%$]{
            \includegraphics[width=0.23\textwidth]{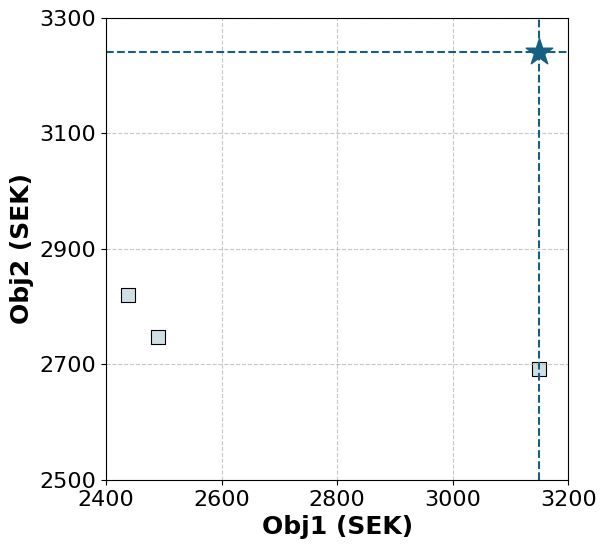}}
            \hspace{0.02in}
            \subfigure[B3M2 with $\epsilon=2\%$]{
            \includegraphics[width=0.23\textwidth]{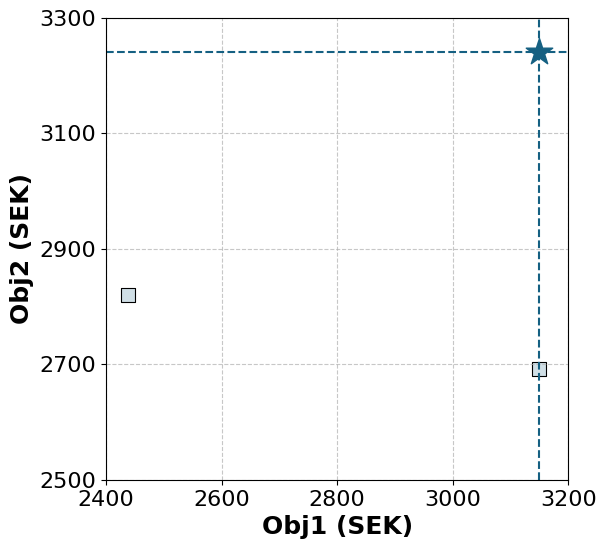}}
            \hspace{0.02in}
            \subfigure[B3M2 with $\epsilon=5\%$]{
            \includegraphics[width=0.23\textwidth]{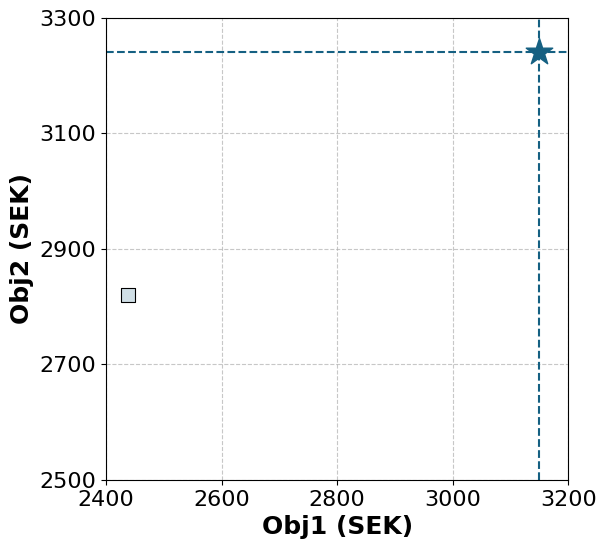}}
        \caption{Non-dominated solutions of Case \#1 with different methods}
	\label{FIG:4}
\end{figure*}

\begin{figure*} [htbp]
	\centering
            \subfigure[All the solutions]{
		  \includegraphics[width=0.23\textwidth]{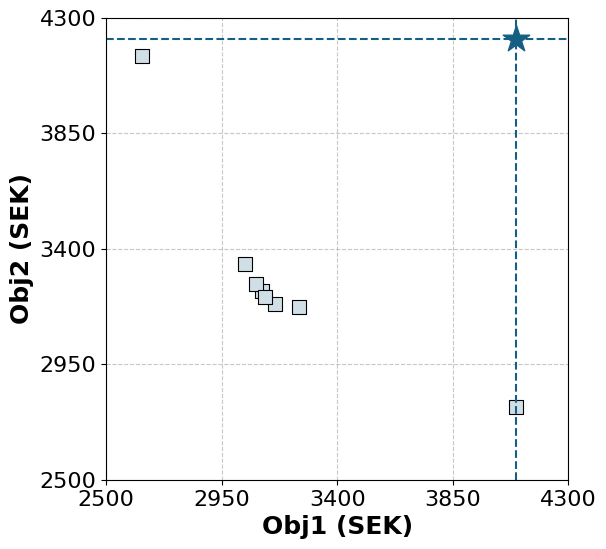}}
            \hspace{0.02in}
            \subfigure[B3M1 with $\epsilon=1\%$]{
            \includegraphics[width=0.23\textwidth]{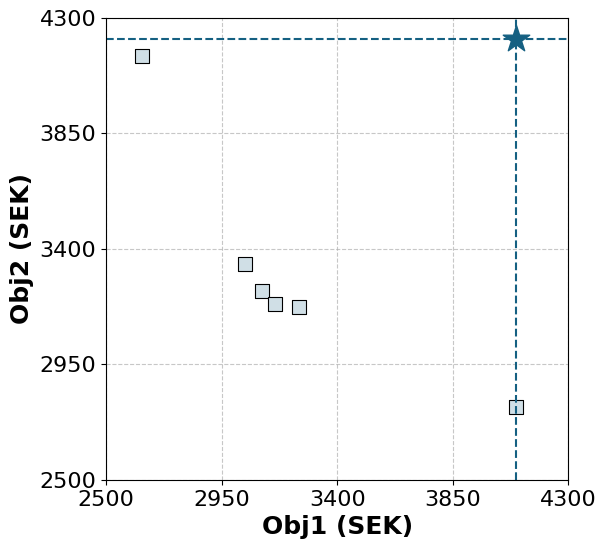}}
            \hspace{0.02in}
            \subfigure[B3M1 with $\epsilon=2\%$]{
            \includegraphics[width=0.23\textwidth]{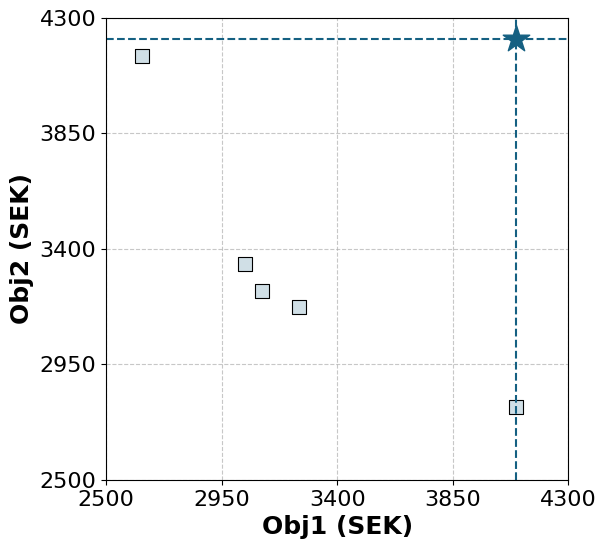}}
            \hspace{0.02in}
            \subfigure[B3M1 with $\epsilon=5\%$]{
            \includegraphics[width=0.23\textwidth]{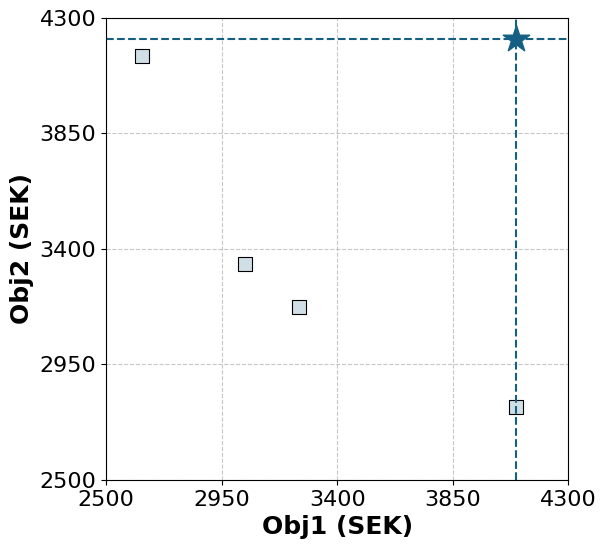}}
            \hspace{0.02in}
            \subfigure[B3M2 with $\epsilon=1\%$]{
            \includegraphics[width=0.23\textwidth]{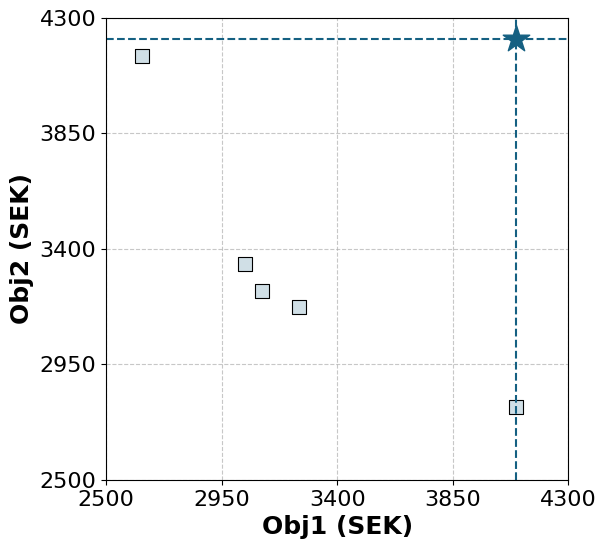}}
            \hspace{0.02in}
            \subfigure[B3M2 with $\epsilon=2\%$]{
            \includegraphics[width=0.23\textwidth]{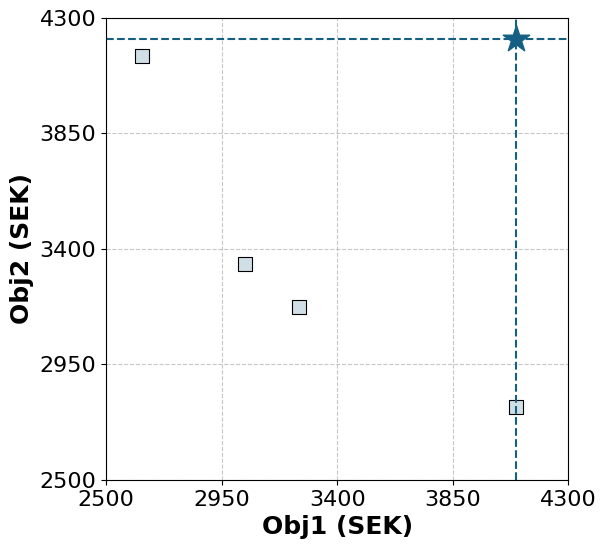}}
            \hspace{0.02in}
            \subfigure[B3M2 with $\epsilon=5\%$]{
            \includegraphics[width=0.23\textwidth]{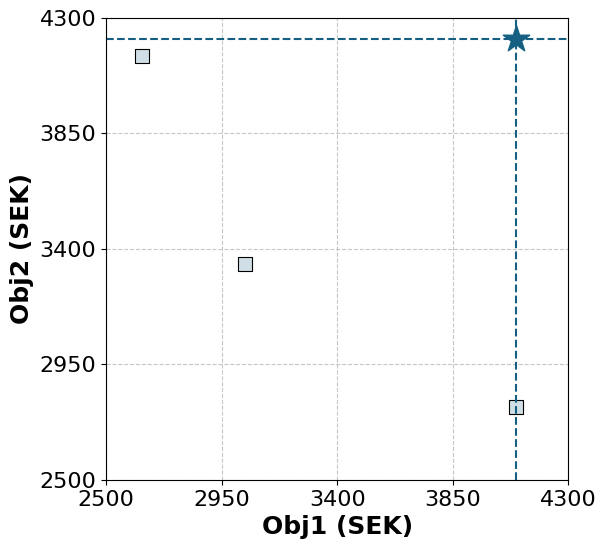}}
        \caption{Non-dominated solutions of Case \#2 with different methods}
	\label{FIG:5}
\end{figure*}

\begin{figure*} [htbp]
	\centering
            \subfigure[All the solutions]{
		  \includegraphics[width=0.23\textwidth]{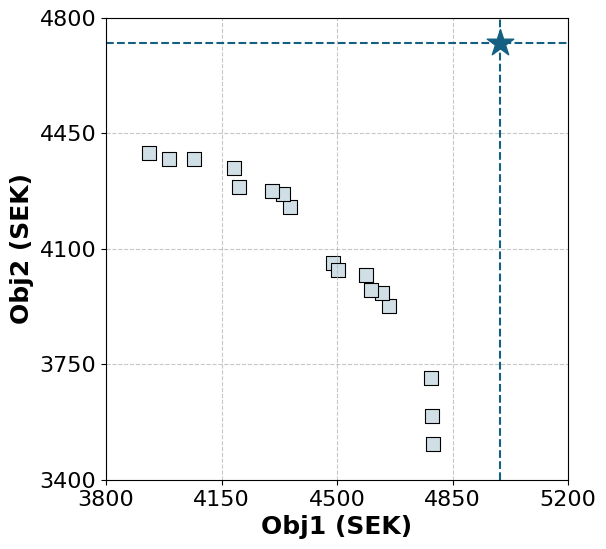}}
            \hspace{0.02in}
            \subfigure[B3M1 with $\epsilon=1\%$]{
            \includegraphics[width=0.23\textwidth]{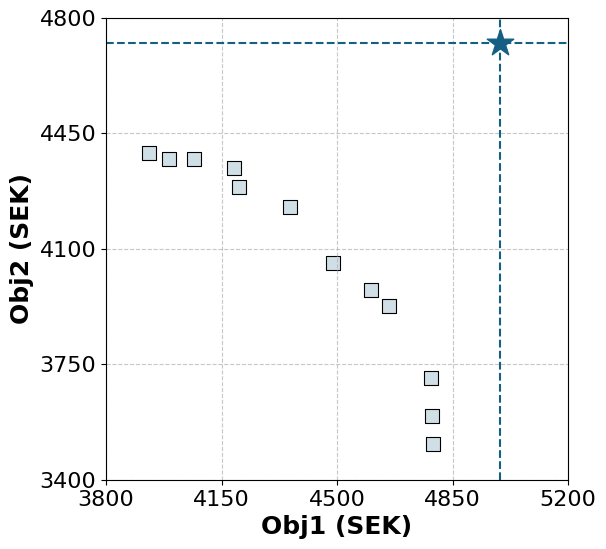}}
            \hspace{0.02in}
            \subfigure[B3M1 with $\epsilon=2\%$]{
            \includegraphics[width=0.23\textwidth]{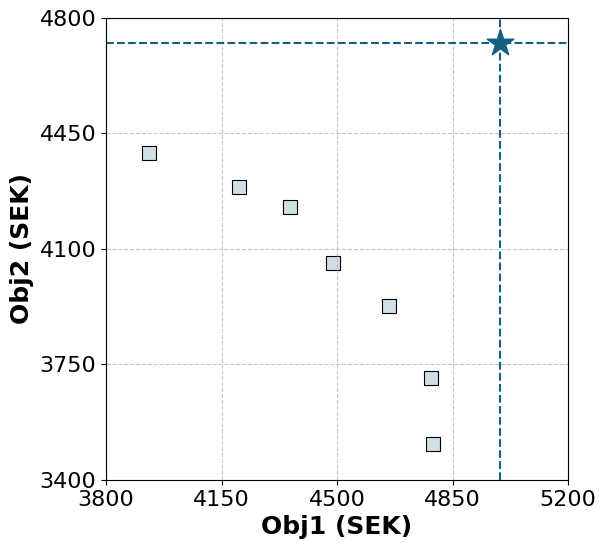}}
            \hspace{0.02in}
            \subfigure[B3M1 with $\epsilon=5\%$]{
            \includegraphics[width=0.23\textwidth]{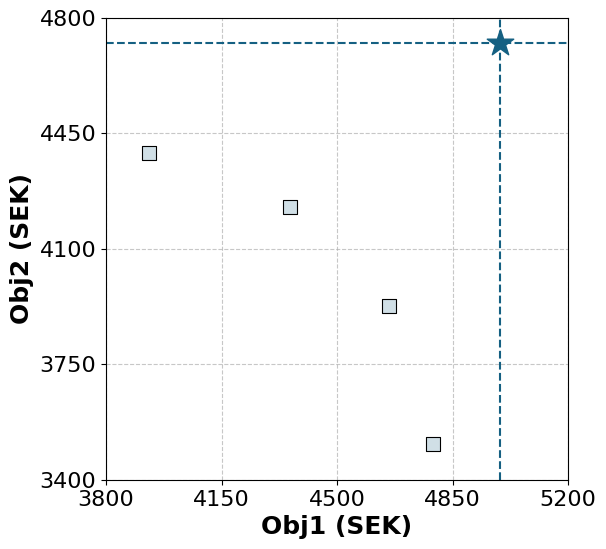}}
            \hspace{2.22in}
            \subfigure[B3M2 with $\epsilon=1\%$]{
            \includegraphics[width=0.23\textwidth]{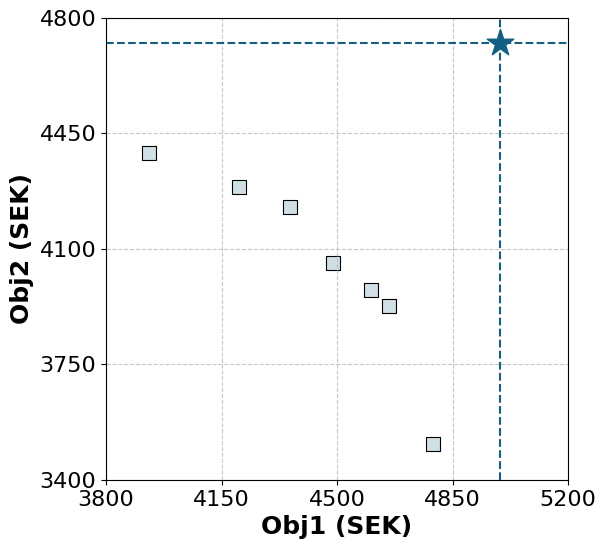}}
            \hspace{0.02in}
            \subfigure[B3M2 with $\epsilon=2\%$]{
            \includegraphics[width=0.23\textwidth]{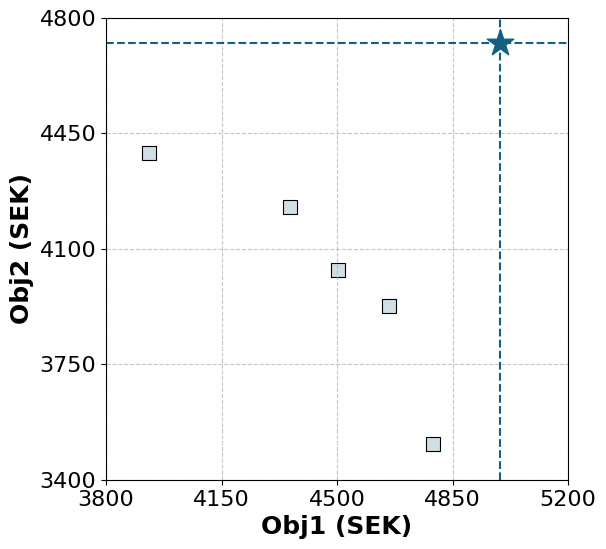}}
            \hspace{0.02in}
            \subfigure[B3M2 with $\epsilon=5\%$]{
            \includegraphics[width=0.23\textwidth]{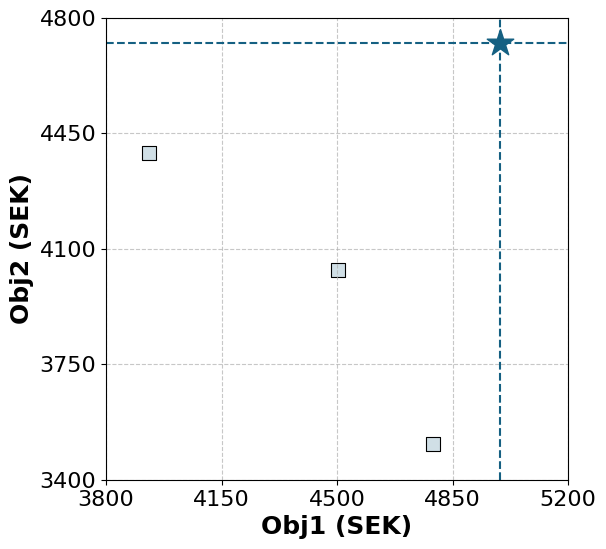}}
        \caption{Non-dominated solutions of Case \#3 with different methods}
	\label{FIG:6}
\end{figure*}

The non-dominated solutions for the three cases are presented in Figs. \ref{FIG:4}, \ref{FIG:5}, and \ref{FIG:6}. The figures illustrate the solutions obtained through the implementation of the original balanced box method and two B3Ms.

As observed in Fig. \ref{FIG:4}, the non-dominated solutions are primarily concentrated in the lower-left region in Case \#1. With B3M1, the non-dominated solutions are distributed in a relatively uniform manner. In contrast, the use of B3M2 results in a reduction in the number of non-dominated points, particularly in scenarios where $\epsilon=2\%$ and $\epsilon=5\%$. This is due to the close proximity of the upper-left and lower-right points in terms of objective 2, limiting the diversity of reduced solutions. It can be concluded that if the variation in one of the objectives is insufficient, then B3M2 becomes ineffective. In such cases, B3M1 is recommended as a more suitable alternative.

In Fig. \ref{FIG:5}, the non-dominated points are mainly concentrated in the central region, forming a compact cluster of solutions. The non-dominated solutions are distributed relatively evenly among the two B3Ms. Given the similarity in the ranges of the two objectives, both methods yielded comparable results, with B3M2 producing a slight reduction in the number of non-dominated points.

As illustrated in Fig. \ref{FIG:6}, the number of non-dominated points is considerable, and no distinct clustering is evident. The distribution of non-dominated points is more uniform with B3M2, as illustrated in the two pairs of comparisons: Fig. \ref{FIG:6} (c) vs. (f), and Fig. \ref{FIG:6} (d) vs. (g). This is due to the fact that B3M2 is more aggressive than B3M1. In contrast to B3M1, which ignores a current solution only when the differences between both objectives and an existing solution fall within the tolerance range, B3M2 is more rigorous in maintaining diversity. This is because B3M2 ignores a solution if either objective falls within the tolerance range. The advantage of B3M2 is that it facilitates a more uniform distribution of non-dominated points, thereby ensuring a greater degree of diversity and coverage across the objective space.

\subsection{Sensitivity analysis}

\begin{figure*} [htbp]
	\centering
            \subfigure[All the solutions]{
            \includegraphics[scale=.27]{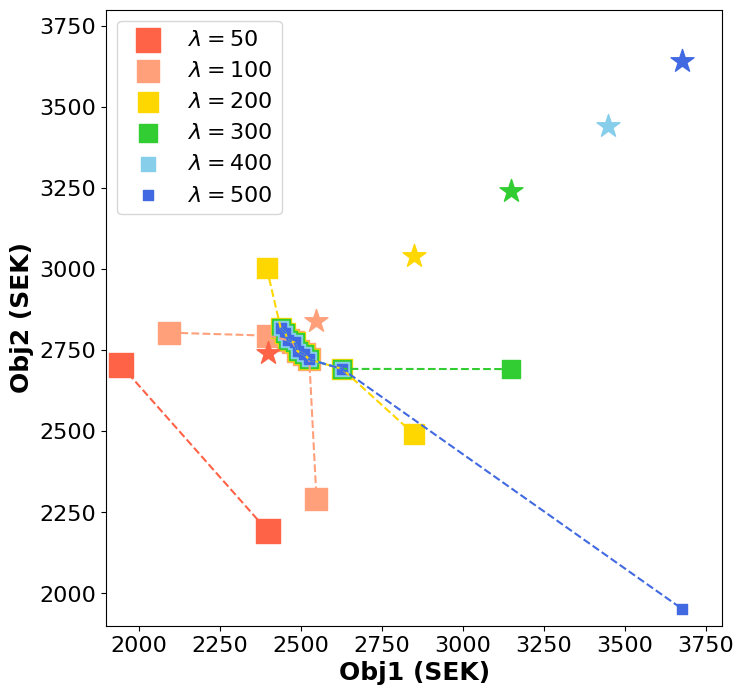}}
            \hspace{0.02in}
            \subfigure[B3M1 with $\epsilon=1\%$]{
            \includegraphics[scale=.27]{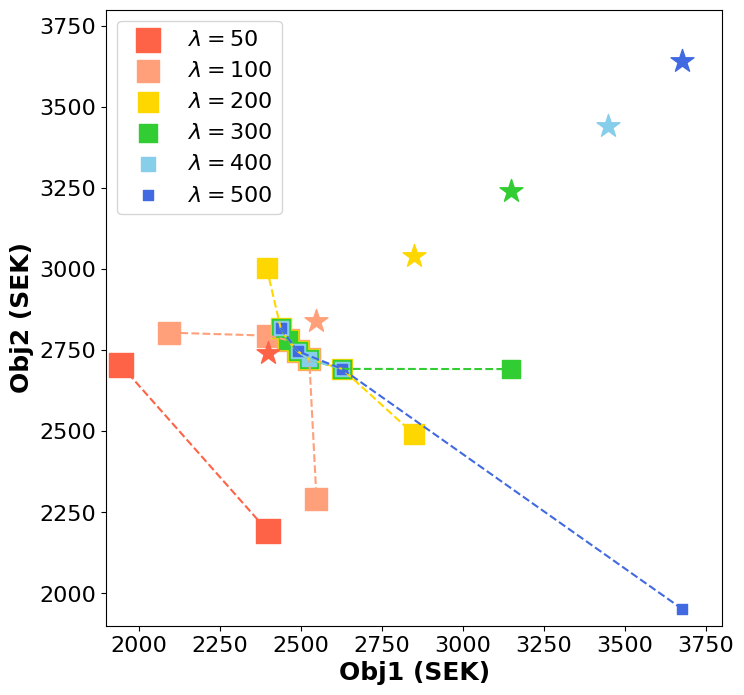}}
            \hspace{0.02in}
            \subfigure[B3M2 with $\epsilon=1\%$]{
            \includegraphics[scale=.27]{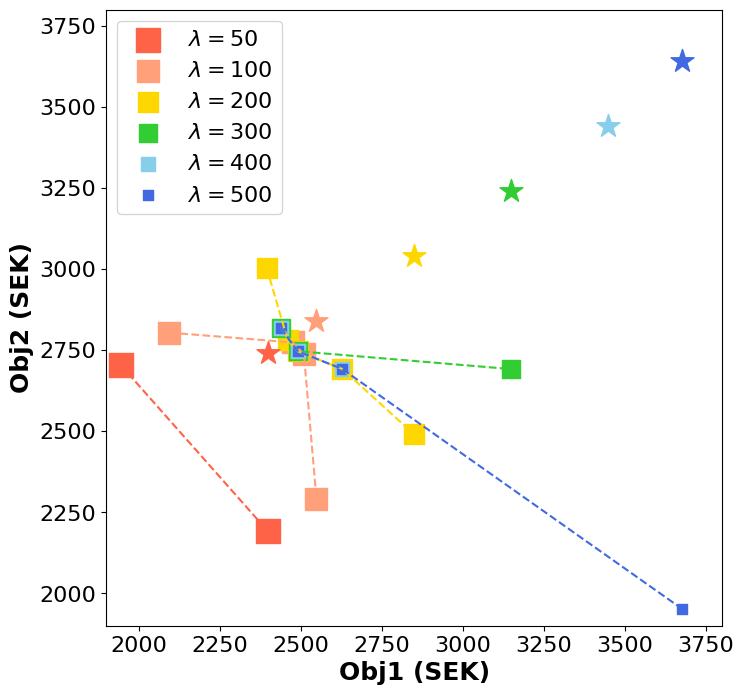}}
        \caption{Impact of VOT levels ($\lambda$) on optimized results using B3Ms}
	\label{FIG:VOTM1}
\end{figure*}

This section examines the effects of varying levels of value of time (VOT, $\lambda$) on the optimization process and the balance between competing objectives. Fig. \ref{FIG:VOTM1} illustrates the comparative performance of collaborative and non-collaborative solutions at varying VOT levels, thereby demonstrating the influence on the optimization outcomes. The squares of different colors represent collaborative solutions; the dashed lines connect them to form a non-dominated frontier; the stars indicate the corresponding non-collaborative solutions.

As shown in Fig. \ref{FIG:VOTM1} (a), an increase in VOT generally results in a more pronounced advantage of collaborative solutions over non-collaborative ones, particularly at higher VOT levels (e.g., $\lambda=500$ SEK). At these levels, non-collaborative solutions deviate significantly from the frontier. This highlights the vital importance of collaboration in achieving superior outcomes under conditions of high VOT.
As VOT increases, the trade-offs between objectives become more intricate, and the frontier becomes steeper, reflecting the growing challenge of decision-making.
Notably, several points overlap across different VOT levels, indicating that certain solutions remain robust or near-optimal despite changes in VOT. These solutions could indicate regions of the solution space where trade-offs are less sensitive to VOT changes, supporting more stable decision-making strategies.
Figs. \ref{FIG:VOTM1} (b) and (c) compare the results of the two B3Ms at different VOT levels. With $\epsilon=1\%$ in B3M1 (Fig. \ref{FIG:VOTM1} (b)), the solutions are more precise, and the distance between collaborative and non-collaborative solutions is greater, especially at higher $\lambda$ values. In contrast, B3M2 (Fig. \ref{FIG:VOTM1} (c)) produces more dispersed solutions, reflecting its more aggressive nature. However, this can compromise precision, particularly at higher $\lambda$ values.
Due to page limitations, we present only the results of B3Ms with $\epsilon=1\%$. As $\epsilon$ increases, the solution set becomes more dispersed, the efficient frontier flattens, and the advantage of collaboration diminishes. Higher $\epsilon$ values also reduce solution accuracy while shortening computation time, highlighting a trade-off between precision and computational efficiency.

\begin{figure*} [htbp]
	\centering
            \subfigure[All the solutions]{
            \includegraphics[scale=.27]{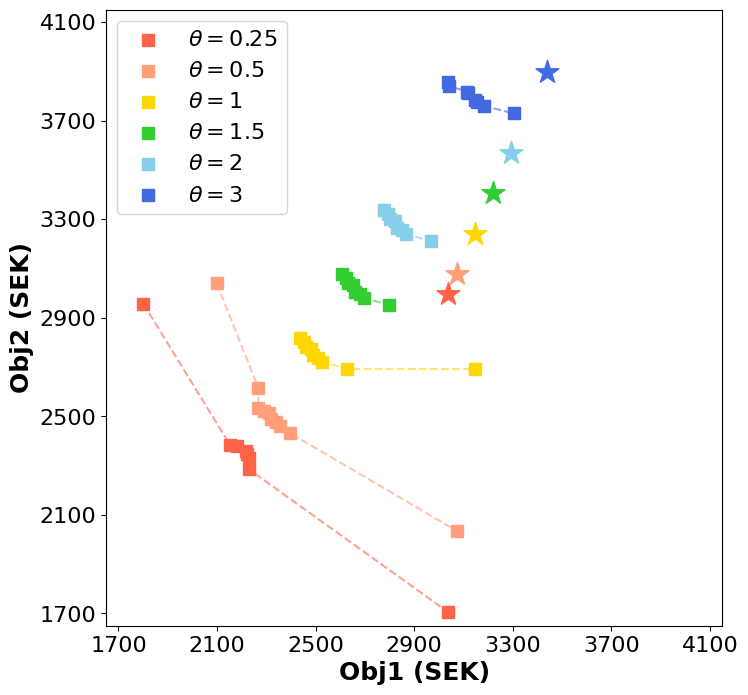}}
            \hspace{0.02in}
            \subfigure[B3M1 with $\epsilon=1\%$]{
            \includegraphics[scale=.27]{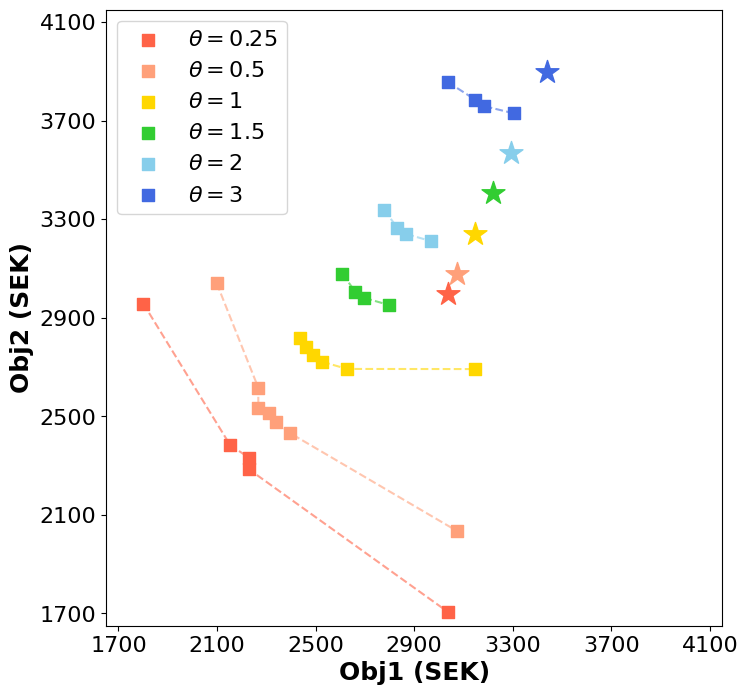}}
            \hspace{0.02in}
            \subfigure[B3M2 with $\epsilon=1\%$]{
            \includegraphics[scale=.27]{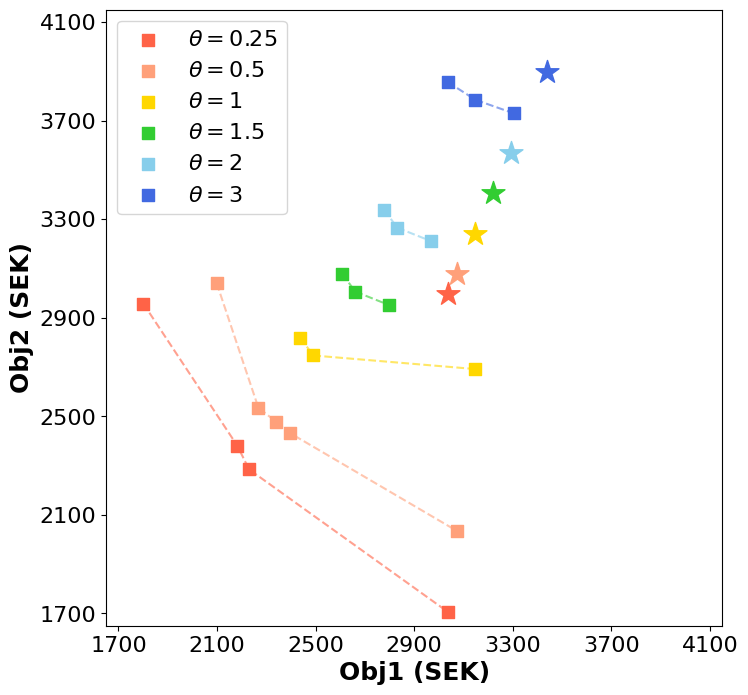}}
        \caption{Impact of changes in unit energy cost ($\theta$) on optimized results using B3Ms}
	\label{FIG:unitECost1}
\end{figure*}

In addition, the energy fees $c_{\rm e,c}^{j,\pi}$ and $c_{\rm e,o}^{j,\pi}$ are adjusted in parallel by a uniform scaling factor $\theta$. The adjusted fees are then calculated as $c=c \times \theta$. This approach allows an investigation of the influence of proportional changes in both collaborative and own EV charging fees on the optimization results. Fig. \ref{FIG:unitECost1} shows the impact of changes in unit energy cost ($\theta$) on both collaborative and non-collaborative solutions.

In Fig. \ref{FIG:unitECost1} (a), at lower unit energy costs (smaller $\theta$), the distance between cooperative and non-collaborative solutions is greater, thereby emphasizing the necessity and advantages of collaboration. Conversely, as unit energy costs increase (larger $\theta$), this distance decreases, thereby diminishing the benefits and importance of collaboration. One potential explanation for this reduction in disparity is that rental costs for charging stations remain unchanged, which may prompt companies to prefer leasing stations on their own as energy costs rise. This may, in turn, result in a reduction in the necessity for collaboration. Furthermore, as the value of $\theta$ increases, the efficient frontier shortens, and the collaborative solutions become more concentrated. This indicates that higher energy costs limit the range of trade-offs between objectives, thereby reducing both the flexibility in optimization and the variety of available solutions.
Figs. \ref{FIG:unitECost1} (b) and (c) further extend the analysis by illustrating the impact of different scaling factors ($\theta$) on the results of B3M1 and B3M2. A similar trend is observed in Fig. \ref{FIG:VOTM1}, with B3M2 continuing to demonstrate a more dispersed set of solutions across the $\theta$ values, reflecting its more aggressive nature. In contrast, B3M1 maintains a more concentrated set of solutions. This comparison serves to reinforce the trade-off between precision and solution flexibility across different B3Ms and error tolerances.

\subsection{Results of Nash bargaining}

\begin{figure*} [htbp]
	\centering
            \subfigure[Generalized Nash bargaining on all solutions]{
		  \includegraphics[width=0.23\textwidth]{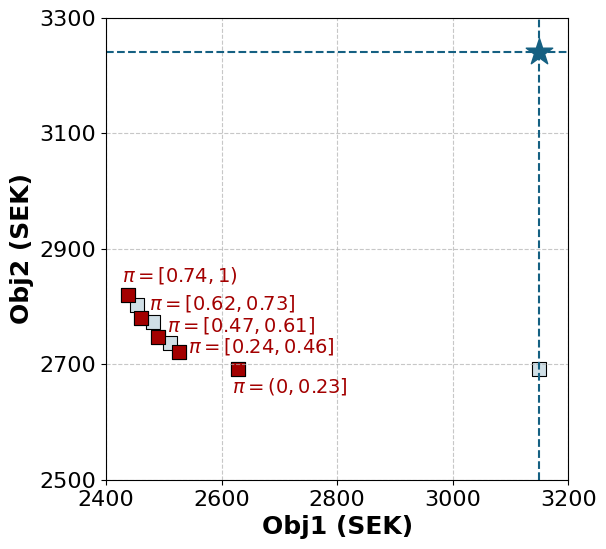}}
            \hspace{0.02in}
            \subfigure[Distance function on all solutions]{
            \includegraphics[width=0.23\textwidth]{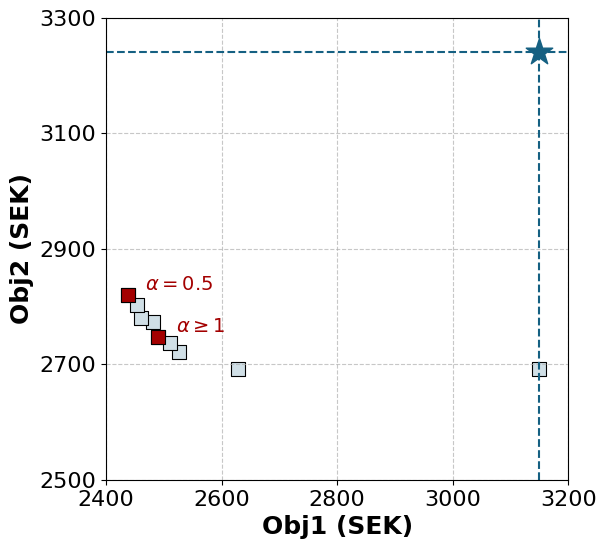}}
            \hspace{0.02in}
            \subfigure[Generalized Nash bargaining on B3M1 with $\epsilon=2\%$]{
		  \includegraphics[width=0.23\textwidth]{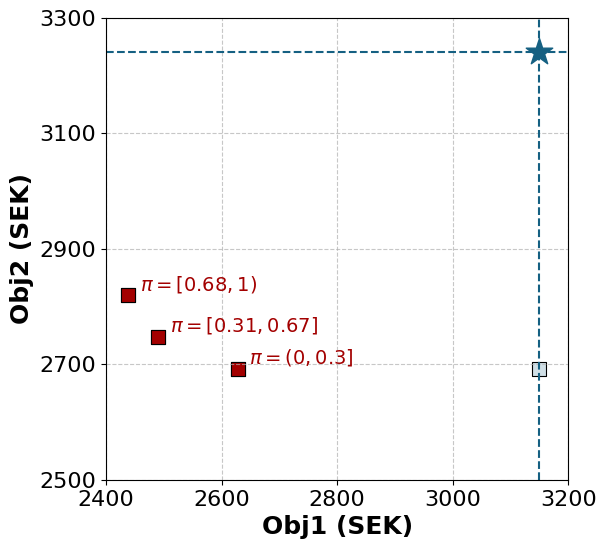}}
            \hspace{0.02in}
            \subfigure[Distance function on B3M1 with $\epsilon=2\%$]{
            \includegraphics[width=0.23\textwidth]{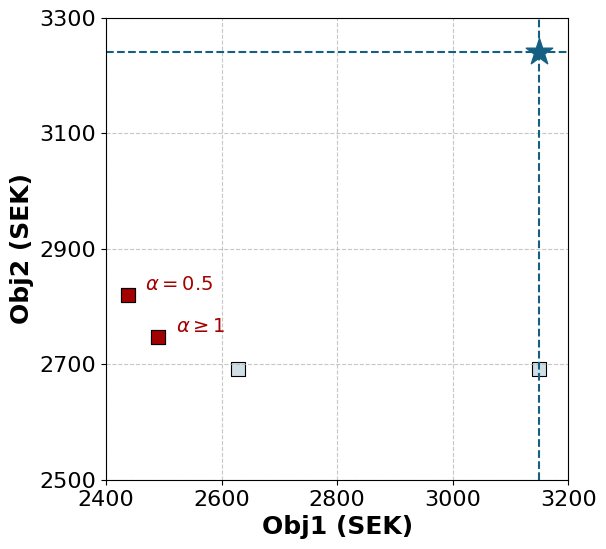}}
        \caption{Final solution of Case \#1 with different methods}
	\label{FIG:8}
\end{figure*}

\begin{figure*} [htbp]
	\centering
            \subfigure[Generalized Nash bargaining on all solutions]{
		  \includegraphics[width=0.23\textwidth]{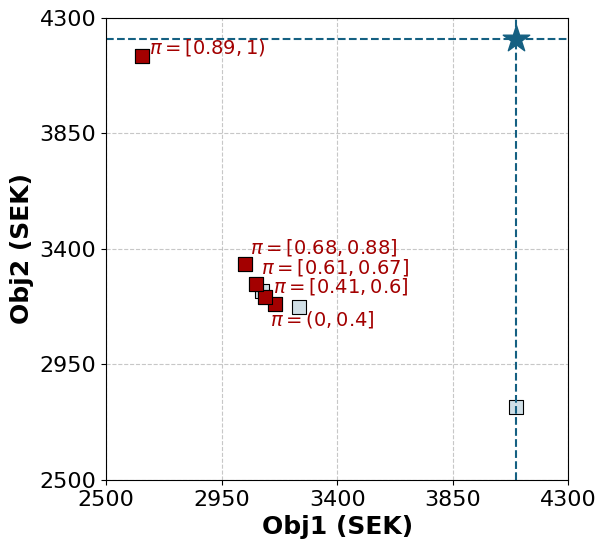}}
            \hspace{0.02in}
            \subfigure[Distance function on all solutions]{
            \includegraphics[width=0.23\textwidth]{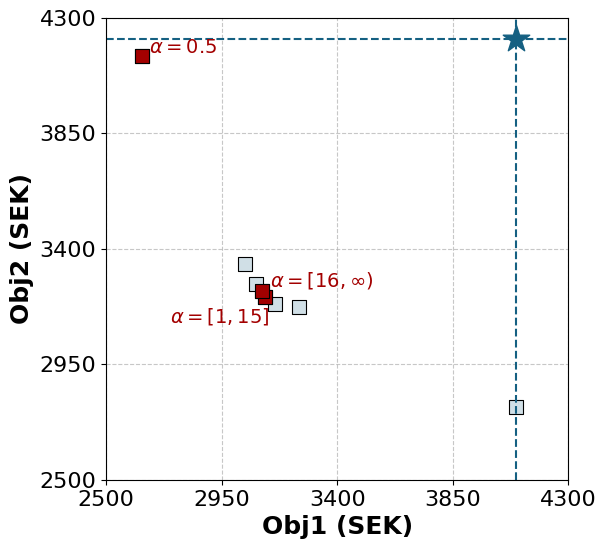}}
            \hspace{0.02in}
            \subfigure[Generalized Nash bargaining on B3M1 with $\epsilon=2\%$]{
		  \includegraphics[width=0.23\textwidth]{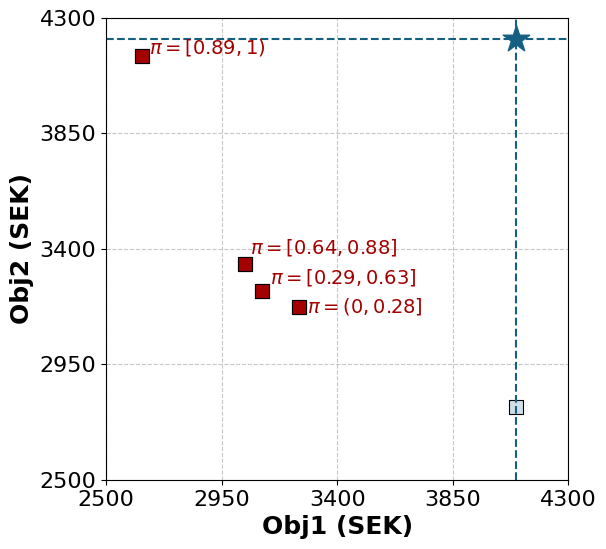}}
            \hspace{0.02in}
            \subfigure[Distance function on B3M1 with $\epsilon=2\%$]{
            \includegraphics[width=0.23\textwidth]{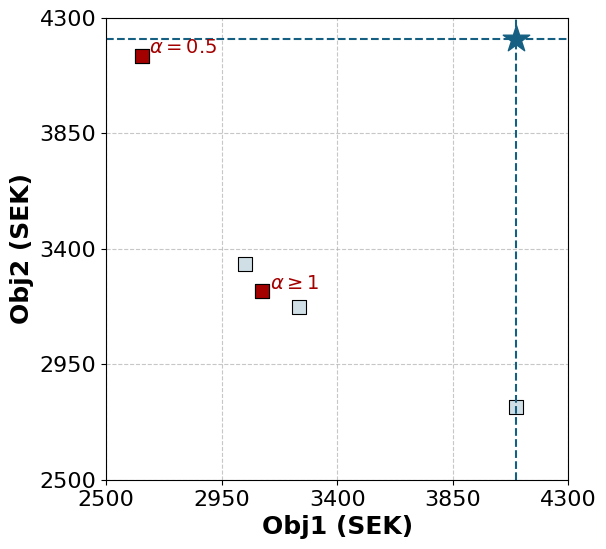}}
        \caption{Final solution of Case \#2 with different methods}
	\label{FIG:9}
\end{figure*}

\begin{figure*} [htbp]
	\centering
            \subfigure[Generalized Nash bargaining on all solutions]{
		  \includegraphics[width=0.23\textwidth]{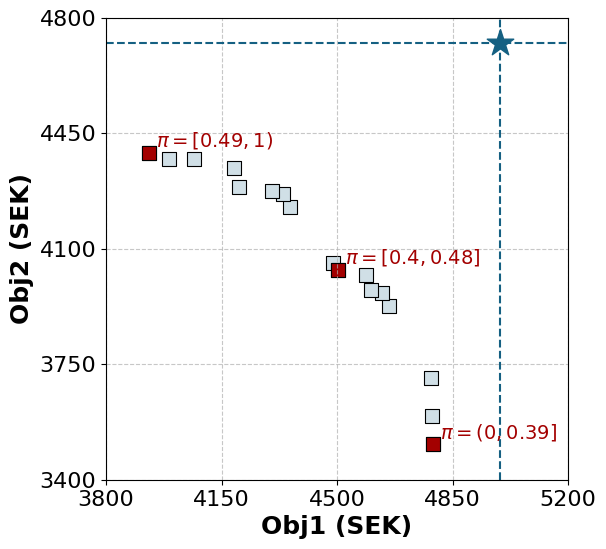}}
            \hspace{0.02in}
            \subfigure[Distance function on all solutions]{
            \includegraphics[width=0.23\textwidth]{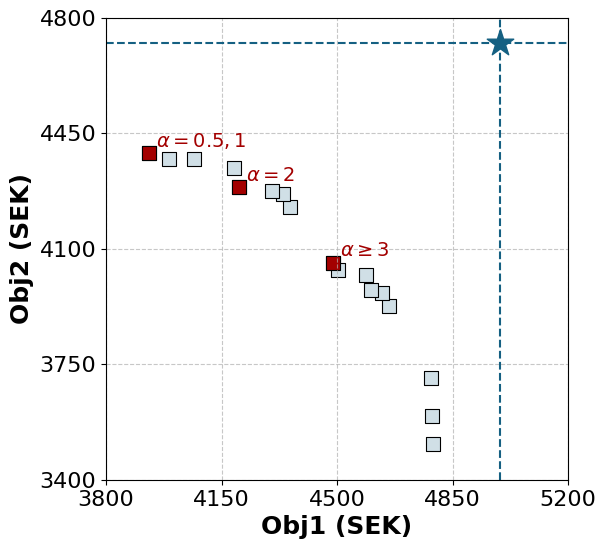}}
            \hspace{0.02in}
            \subfigure[Generalized Nash bargaining on B3M1 with $\epsilon=2\%$]{
		  \includegraphics[width=0.23\textwidth]{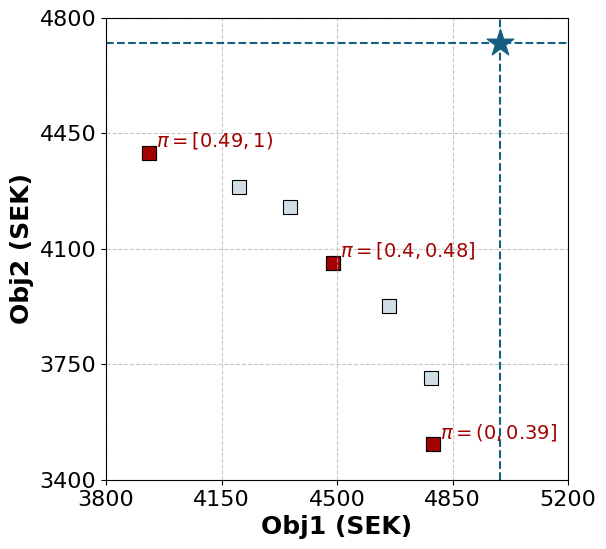}}
            \hspace{0.02in}
            \subfigure[Distance function on B3M1 with $\epsilon=2\%$]{
            \includegraphics[width=0.23\textwidth]{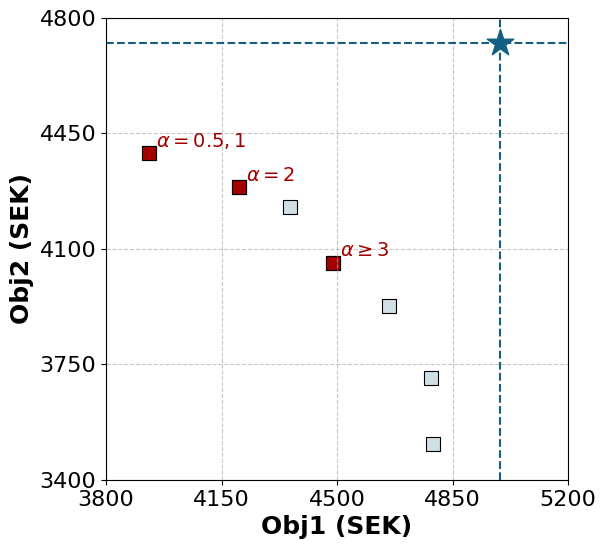}}
        \caption{Final solution of Case \#3 with different methods}
	\label{FIG:10}
\end{figure*}

The final solution is identified through the utilization of the generalized Nash bargaining and distance function. The final solutions are illustrated in Figs. \ref{FIG:8}, \ref{FIG:9}, and \ref{FIG:10}, respectively, for different values of the parameters. In the figure, the square markers represent non-dominated points in the solution space. The red squares indicate final solutions for specific parameter intervals, with the parameter values labeled next to each point. These red squares demonstrate the impact of parameter ranges on the trade-offs between objectives. A comparison of the final solutions based on the exact efficient frontier and the efficient frontier indicates no significant discrepancy between them, particularly when the distance function is utilized to identify the final agreement point.

Regarding the generalized Nash bargaining method, it is evident that when the power value $\pi$ is relatively low, player 2 makes a significantly larger claim. Consequently, player 2 incurs a lower cost, whereas player 1 faces a higher cost. 
As the power $\pi$ increases, the final agreement point shifts towards the upper left.
The Nash bargaining solution leverages the parties' relative bargaining power to identify an optimal agreement. By adjusting the relative bargaining power, decision-makers can strategically influence the final agreement, thereby ensuring that it reflects the balance of interests and power dynamics between the parties.
Additionally, it was observed that the efficient frontier diminished the sensitivity of power, a trade-off that enhances computational efficiency but may impact the precision of the Nash bargaining outcome.

Regarding the distance function method, when the sensitivity of the norm $\alpha$ is relatively small (e.g., $\alpha=1$), the errors associated with each objective are given equal weighting, reflecting a balanced approach. Conversely, when $\alpha$ is relatively large (e.g., $\alpha \to \infty$), the optimization focuses on minimizing the largest normalized error, emphasizing robustness by prioritizing the worst-case scenario. In this case, the final agreement point represents the robust optimal solution, ensuring neither objective is excessively compromised. Notably, as shown in Fig. \ref{FIG:8} (b) and (d), the solutions for $\alpha=1$ and $\alpha \to \infty$ yield the same final solution, demonstrating the balance achieved between minimizing total error and addressing the worst deviation. The flexibility of $\alpha$ enables decision-makers to balance objectives and tailor strategies, making it a powerful tool for managing trade-offs in bi-objective optimization.

Both the Nash bargaining method and the distance function ($\alpha$-norm) approach provide frameworks for optimizing bi-objective problems, each offering unique flexibility. The Nash bargaining method considers the relative bargaining power of each party, optimizing outcomes based on their strategic influence. In contrast, the distance function ($\alpha$-norm) allows for the balancing of objectives by adjusting the trade-off between total error and maximum deviation.

\subsection{Charging scheduling results}

\begin{figure*} [htbp]
	\centering
		\includegraphics[scale=.15]{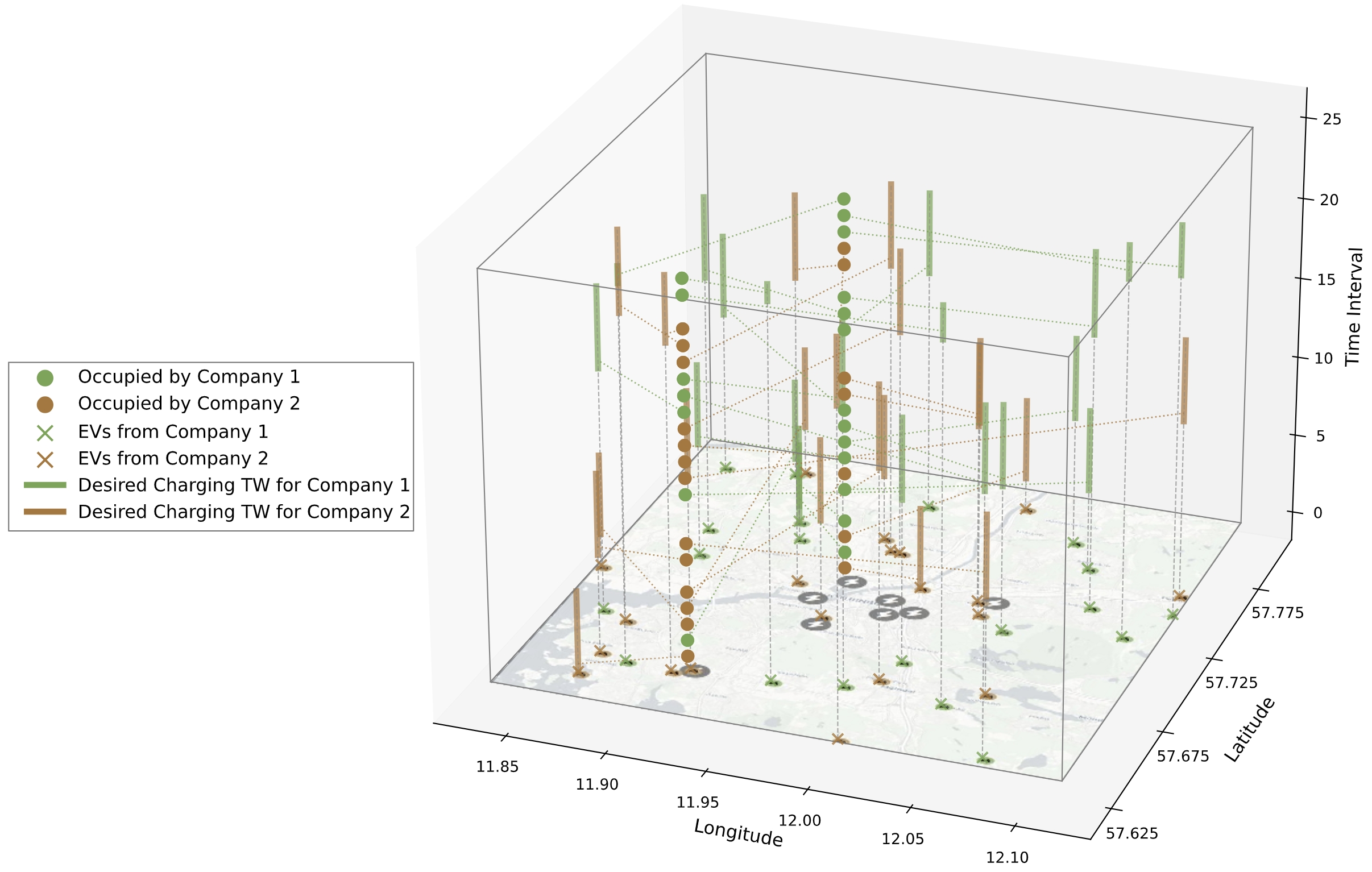}
    \caption{Spatiotemporal distribution of EV charging scheduling results}
	\label{FIG:Spatiotemporal}
\end{figure*}

Figure \ref{FIG:Spatiotemporal} illustrates the spatiotemporal distribution of the optimized charging schedules for Case \#3 ($\pi=\left[0.4-0.48\right]$ or $\alpha \geq 3$) using B3M1.
The results demonstrate how EVs from each company are allocated to specific charging stations during defined time intervals, ensuring that the charging sessions fall within their designated charging TWs. The colors green and orange represent Company 1 and Company 2, respectively, while the "X" markers indicate the positions of the EVs.

As shown in Fig. \ref{FIG:Spatiotemporal}, the optimization results indicate that only two charging stations are rented, a direct result of the high rental costs. EVs from Company 2 are predominantly served by the left-side charging station, while those from Company 1 primarily use the right-side charging station. This preference is attributable to the fact that charging at their own rented stations entails lower electricity costs than using other stations, which is in accordance with the cost-minimization objective of the model.
Notably, the charging window - TW coupling reveals that most EVs begin charging at the start of their respective TWs. This behavior is driven by the incorporation of waiting time costs into the objective function, which motivates the results to minimize delays. The coupling shows that the model works well at making scheduling decisions that are in line with time constraints, which increases its usefulness in real life.
The results show that the proposed model efficiently allocates charging resources and resolves potential conflicts. The absence of overlapping charging assignments and optimal utilization of charging stations demonstrates this.
These findings highlight the robustness of the optimization strategy in addressing the charging demands of both companies, offering valuable insights for enhancing real-world charging management systems.

\section{Conclusion} \label{Section5}

This paper addresses the collaborative scheduling problem in which two companies share charging stations, utilizing a bi-objective optimization approach to minimize costs for both parties through a collaborative scheduling model. We present the balanced box method to generate all non-dominated solutions and develop B3Ms to efficiently produce an efficient frontier with lower computational effort. Furthermore, the Nash bargaining approach is used to ensure equitable cooperation between the companies. Overall, the proposed bi-objective model and developed methods demonstrate their versatility and broad applicability, contributing to the advancement of sustainable urban logistics through horizontal collaboration.

The B3Ms proposed in this study can be applied to all bi-objective integer programming models. These methods are promising because they can significantly reduce computational time while preserving the structure of the efficient frontier. We introduced two distinct methods: the first is more conservative, involving searching for a solution within a defined rectangle and disregarding those that are very similar to existing solutions. In contrast, the second method is more aggressive, directly shrinking the rectangular space to eliminate certain solutions. While the aggressive method may miss some solutions compared to the conservative one, leading to gaps around the predetermined tolerance range, it is more efficient in terms of computational time. The conservative method maintains solutions well within the tolerance range, sometimes even tighter. In addition, the cooperative bargaining method we propose can be broadly applied to all bi-objective programming models to identify the final agreement point. Although the two approaches differ in their proximity to key reference points, they ultimately converge toward the same direction: one method is positioned farther from the non-dominated point, while the other is closer to the ideal point. The two companies can set parameter values, such as power, and negotiate to reach the final solution.

Beyond the core contributions of this study, several promising directions remain open for further research. First, scaling the method to very large problem instances remains a practical challenge, as even a reduced set of exact representative solutions may become computationally intensive. Developing region-wise decomposition strategies or distributed versions of B3Ms could improve tractability while preserving solution interpretability and structural properties. Second, incorporating real-time operational uncertainties, such as unpredictable vehicle arrivals, dynamic fleet participation, and fluctuations in station availability, would improve the robustness and adaptability of the scheduling framework under realistic conditions. Finally, extending the model to accommodate interruptive scheduling behaviors, where EVs may be preempted, reassigned, or rescheduled in response to operational disruptions or shifting priorities, would further increase its flexibility and practical relevance in complex fleet operations. These directions collectively aim to improve the scalability, adaptability, and decision-support capability of collaborative EV scheduling models in increasingly dynamic transportation environments.

\section*{Acknowledgments}
This work was supported by the European Commission, Swedish Energy Agency, and Chalmers University of Technology through the projects E-Laas (Energy optimal urban Logistics As A Service) and COLLECT (Horizontal Cooperation in Urban Distribution Logistics - a Trusted - Cooperative Electric Vehicle Routing Method). Partially supported by the Swedish Electromobility Center (SEC) under the project LEAR: Robust LEArning methods for electric vehicle Route selection.

\appendix

\section*{Appendix} \label{Appendix}

\setcounter{figure}{0}  
\renewcommand{\thefigure}{A\arabic{figure}}

\begin{figure*} [htbp]
	\centering
		\includegraphics[scale=0.58]{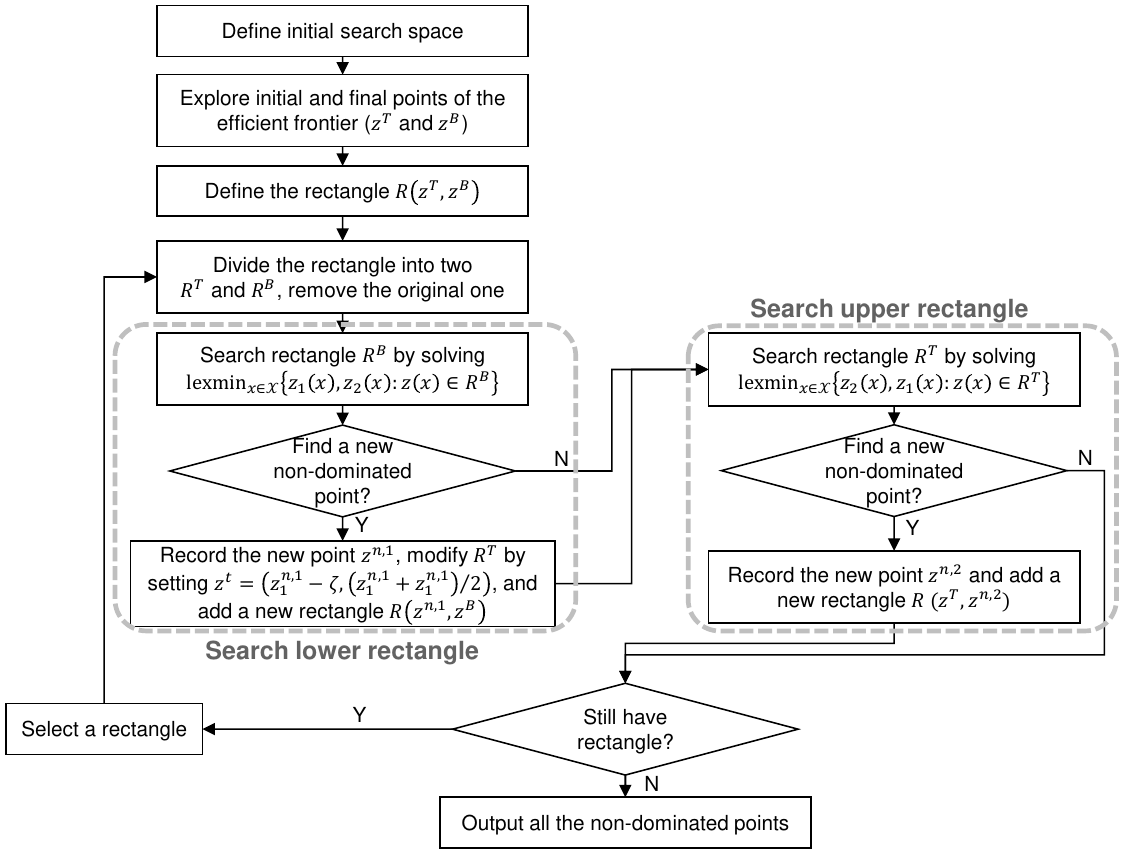}
    \caption{The flowchart of the balanced box method}
	\label{FIG:flowchart}
\end{figure*}

The flowchart for the balanced box method is illustrated in Fig \ref{FIG:flowchart}. This method revolves around the exploration of rectangles and comprises four main steps. A comprehensive explanation of these steps is provided below, offering further clarification on the operations in the flowchart.

\textit{Step 1}: To identify a yet unidentified non-dominated point, it is necessary to divide the rectangle $R(z^T,z^B)$ into two distinct regions, the top rectangle ($R^T$) and the bottom rectangle ($R^B$). These two rectangles are defined by the points $z^T$, $z^t$, $z^b$, and $z^B$. In this case, the value of $z^t$ is $(z^B_1, (z^T_2+z^B_2)/2)$, while $z^b$ is $(z^T_1, (z^T_2+z^B_2)/2)$. In other words, the original rectangle is divided horizontally along the $z_2(\mathscr{x})$ axis, and subsequently removed from the rectangle set.

\textit{Step 2}: The lower rectangle, $R^B(z^b,z^B)$, can be searched by solving the following optimization problem: $\mathop{\rm{lexmin}}_{\mathscr{x}\in \mathscr{X}} \left\{z_1(\mathscr{x}),z_2(\mathscr{x}):z(\mathscr{x})\in R^B\right\}$. If a new non-dominated point is found, then a portion of the rectangle $R^T(z^T,z^t)$ is dominated. The rectangle $R^T(z^T,z^t)$ should be modified by setting $z^t$ equal to $(z^{\rm n,1}_1-\zeta,(z^T_2+z_2^B)/2)$, where $\zeta$ is a small constant. Furthermore, a new rectangle, $R(z^{\rm n,1},z^B)$, should be incorporated into the existing rectangle set. The non-dominated point is then recorded. If no new solution is identified, $R^T(z^T,z^t)$ remains unchanged, and no further rectangle is incorporated.

\textit{Step 3}: The upper rectangle, $R^T(z^T,z^t)$, can be searched by solving the following optimization problem: $\mathop{\rm{lexmin}}_{\mathscr{x}\in \mathscr{X}} \left\{z_2(\mathscr{x}),z_1(\mathscr{x}):z(\mathscr{x})\in R^T\right\}$. If a new non-dominated point is identified, then a new rectangle, $R(z^T,z^{\rm n,2})$, should be incorporated into the existing rectangle set. The non-dominated point is then recorded; otherwise, no new rectangle is included.

\textit{Step 4}: If the rectangle set is not empty, a single rectangle should be selected from the set and proceed with step 1; otherwise, all the non-dominated points must be outputted.

\bibliographystyle{cas-model2-names}

\bibliography{reference}

\end{document}